\documentclass[reqno,11pt,letterpaper]{amsart}
\pdfoutput=1
\numberwithin{equation}{section}

\usepackage[utf8]{inputenc}
\usepackage[british]{babel}

\usepackage{amsmath}
\usepackage{amsfonts}
\usepackage{amssymb,esint}
\usepackage{tikz}
\usepackage{bigints}
\usepackage{comment}
\usepackage{geometry}
\usepackage{lbmath}
\usepackage[autostyle]{csquotes}

\usepackage{nccmath}
\usepackage{color}
\usepackage{mathrsfs}
\usepackage{cancel}
\usepackage{mathtools}
\usepackage{abs}
\usepackage{set}
\usepackage{norm}
\usepackage{seminorm}
\usepackage[angle-vert]{ip}
\usepackage{calligra}

\usepackage{etoolbox}
\cslet{blx@noerroretextools}\empty
\usepackage[maxbibnames=9,giveninits=true]{biblatex}
\DeclareFieldFormat*{title}{#1}
\usepackage[colorlinks, breaklinks, citecolor=citegreen, linkcolor=refred,urlcolor=blue, unicode,psdextra]{hyperref}

\usepackage{bookmark}
\usepackage{cleveref}
\expandafter\def\csname ver@etex.sty\endcsname{3000/12/31}

\usepackage{autonum}
\definecolor{mentadent}{RGB}{0,122,41}
\definecolor{citegreen}{rgb}{0,0.3,0}
\definecolor{refred}{rgb}{0.5,0,0}

\theoremstyle{plain}
\newtheorem{theorem}{Theorem}[section]
\newtheorem{lemma}[theorem]{Lemma}
\newtheorem{proposition}[theorem]{Proposition}
\newtheorem{corollary}[theorem]{Corollary}
\newtheorem{definition}[theorem]{Definition}
\newtheorem{remark}[theorem]{Remark}

\theoremstyle{remark}
\crefname{definition}{Definition}{Definitions}
\crefname{theorem}{Theorem}{Theorems}
\crefname{lemma}{Lemma}{Lemmas}
\crefname{step}{Step}{Steps}
\crefname{substep}{Step}{Steps}
\crefname{proposition}{Proposition}{Propositions}
\crefname{corollary}{Corollary}{Corollaries}
\crefname{remark}{Remark}{Remarks}
\crefname{section}{Section}{Sections}
\crefname{subsection}{Section}{Sections}
\crefname{enumi}{}{}

\newcommand{\R}{\mathbb R}
\newcommand{\N}{\mathbb N}

\newcommand{\Int}{\mathrm{Int}}

\let\otheta\theta
\renewcommand{\theta}{\vartheta}

\newcommand{\CC}{\mathscr{C}}
\newcommand{\Cn}{\mathrm{C}}

\DeclareMathOperator{\Crit}{Crit}

\newcommand{\barint}
{\rule[.036in]{.12in}{.009in}\kern-.16in \displaystyle\int\limits}

\let\div=\relax
\DeclareMathOperator{\div}{div}
\newcommand{\pa}{\partial }

\newcommand{\ee}{\mathrm{e}}
\let\oldchi=\chi
\renewcommand{\chi}{\raisebox{\depth}{\(\oldchi\)}}

\usepackage{color}

\newcommand{\numberset}{\mathbb}
\renewcommand{\N}{\numberset{N}}
\renewcommand{\R}{\numberset{R}}
\newcommand{\Sf}{\numberset{S}}

\newcommand{\loc}{\mathrm{loc}}

\newcommand{\ric}{\mathop {\rm Ric}\nolimits}

\DeclareMathOperator{\AVR}{AVR}
\let\d=\relax
\newcommand{\d}[1][g]{d_{#1}}
\newcommand{\capa}{{\rm Cap}}
\newcommand{\D}{{\rm D}}
\newcommand{\dd}{{\,\rm d}}

\newcommand{\HH}{{\rm H}}
\newcommand{\hh}{{\rm h}}

\newcommand{\Om}{\Omega}
\renewcommand{\phi}{\varphi}
\renewcommand{\epsilon}{\varepsilon}
\newlength{\listalpha}
\renewenvironment{cases}[1][l@{\ \ }l]{\arraycolsep=1.4pt\left\lbrace\kern-3pt\begin{array}{#1}}{\end{array}\right.}

\makeatletter
\newcommand*{\getlength}[1]{\strip@pt#1}

\makeatother
\title[Minkowski Inequality on manifolds with nonnegative Ricci]{Minkowski Inequality on complete Riemannian manifolds with nonnegative Ricci curvature}

\author[L.~Benatti]{Luca Benatti}
\address{L.~Benatti, Universit\`a degli Studi di Trento,
via Sommarive 14, 38123 Povo (TN), Italy}
\email{luca.benatti@unitn.it}

\author[M.~Fogagnolo]{Mattia Fogagnolo}
\address{M.~Fogagnolo, Centro di Ricerca Matematica Ennio De Giorgi, Scuola Normale Superiore,
Piazza dei Cavalieri 3, 56126 Pisa (PI), Italy}
\email{mattia.fogagnolo@sns.it}

\author[L.~Mazzieri]{Lorenzo Mazzieri}
\address{L.~Mazzieri, Universit\`a degli Studi di Trento,
via Sommarive 14, 38123 Povo (TN), Italy}
\email{lorenzo.mazzieri@unitn.it}

\begin{document}

\begin{abstract}
In this paper we consider Riemannian manifolds of dimension at least $3$, with nonnegative Ricci curvature and Euclidean Volume Growth. For every open bounded subset with smooth boundary we establish the validity of an optimal Minkowski Inequality. We also characterise the equality case, provided the domain is strictly outward minimising and strictly mean convex. Along with the proof, we establish in full generality sharp monotonicity formulas, holding along the level sets of $p$-capacitary potentials in $p$-nonparabolic manifolds with nonnegative Ricci curvature.
\end{abstract}
\maketitle

\medskip


\section{Introduction}

\subsection{Statements of the main results} Given an open bounded convex domain with smooth boundary $\Omega\subseteq \R^n$, $n\geq3$, the classical Minkowski Inequality, originally proven in~\cite{Minkowski1903}, gives a sharp lower bound for the average of the mean curvature $\HH$ of $\partial \Omega$ in terms of the inverse of its surface radius, that is
\begin{equation}
\left(\frac{\abs{\S^{n-1}}}{\abs{\partial \Omega}}\right)^{\!\frac{1}{n-1}}\leq \fint\limits_{\partial \Omega} \frac{\HH}{n-1} \dd \sigma,
\end{equation}
with the equality satisfied if and only if $\Omega$ is a ball. It was clear to many authors that such inequality deserved to be further investigated. For example one  would like to relax the convexity assumption on one hand, and to prove that the inequality holds on more general ambient manifolds on the other. 

The first question has been positively answered using techniques based on geometric flows \cite{Huisken}, Optimal Transport \cite{Chang2013,Castillon2010} and recently also Nonlinear Potential Theory \cite{Fogagnolo2019,Agostiniani2019}. The latter  method actually provides the most general statement available so far, namely the Extended Minkowski Inequality, holding for every open bounded domain $\Omega \subseteq \R^n$ with smooth boundary
\begin{equation}\label{eq:ExtendedMinkowskiIntro}
\left(\frac{\abs{\partial \Omega^*}}{\abs{\S^{n-1}}}\right)^{\!\frac{n-2}{n-1}}\leq \frac{1}{\abs{\S^{n-1}}}\int\limits_{\partial \Omega} \abs{\frac{\HH}{n-1}} \dd \sigma \, .
\end{equation}
Here $\Omega^*$ denotes the \emph{strictly outward minimising hull} of $\Omega$. The precise definition of $\Omega^*$ is reported in~\eqref{eq:def_SOMH} below and analysed in full detail in~\cite{Fogagnolo2020a}. However, in this preliminary discussion, we just point out that $ \Omega^*$ minimises the perimeter among bounded subsets containing $\Omega$. 

Many improvements can be found in the literature also concerning the question of extending the Minkowski Inequality to more general settings. Firstly Gallego and Solanes in~\cite{Gallego2005} established quermassintegral inequalities for convex domains in the  Hyperbolic space. Using the Inverse Mean Curvature Flow (IMCF for short), de Lima and Gir\~{a}o in~\cite{Lima2016} extended the result to starshaped and strictly mean-convex domains lying in the same ambient manifold. The IMCF has been also employed to establish a Minkowski-type inequality for outward minimising set sitting in the Schwarzschild manifold by Wei \cite{Wei2018}, in the anti-de Sitter-Schwarzschild manifold by Brendle, Hung and Wang~\cite{Brendle2016}, and on  Asymptotically Flat  Static manifolds by McCormick~\cite{McCormick2017}. 

A natural context where to test the validity of a Minkowski Inequality is provided by complete noncompact Riemannian manifold with nonnegative Ricci curvature. 
A very recent work by Brendle~\cite{Brendle2020} is actually pointing in this direction. Indeed, choosing $f=1$ in~\cite[Corollary 1.5]{Brendle2020} a non sharp Minkowski Inequality can be deduced for complete Riemannian manifolds with nonnegative sectional curvature and Euclidean Volume Growth. In the present paper, we prove the following theorem.

\begin{theorem}[Extended Minkowski Inequality]\label{thm:EMI_maintheorem}
Let $(M,g)$ be complete Riemannian manifold with $\Ric \geq 0$ and Euclidean Volume Growth. Let $\Omega \subseteq M$ be a open bounded set with smooth boundary. Then
\begin{equation}\label{eq:EMI_intro}
\left(\frac{\abs{\partial \Omega^*}}{\abs{\S^{n-1}}}\right)^{\kern-.05cm\frac{n-2}{n-1}}\kern-.1cm\AVR(g)^{\frac{1}{n-1}}\leq \frac{1}{\abs{\S^{n-1}}}\int\limits_{\partial \Omega} \abs{ \frac{\HH}{n-1}} \dd \sigma,
\end{equation}
where $\AVR(g)$ is the asymptotic volume ratio of $(M,g)$, $\HH$ is the mean curvature of $\partial \Omega$ with respect to the outward normal unit vector and $\Omega^*$ is the strictly outward minimising hull of $\Omega$. 
\end{theorem}
In case a strictly outward minimising $\Omega \subset$ with strictly mean-convex boundary achieves the identity in \eqref{eq:EMI_intro}, we show that $M \smallsetminus \Omega$ splits as a (truncated) cone.
\begin{theorem}[Rigidity for the Minkowski Inequality]
\label{thm:rigidity}
A bounded strictly outward minimising $\Omega \subset M$ with smooth strictly mean-convex boundary satisfies
\begin{equation}\label{eq:EMI_intro-identity}
\left(\frac{\abs{\partial \Omega}}{\abs{\S^{n-1}}}\right)^{\kern-.05cm\frac{n-2}{n-1}}\kern-.1cm\AVR(g)^{\frac{1}{n-1}} = \frac{1}{\abs{\S^{n-1}}}\int\limits_{\partial \Omega} { \frac{\HH}{n-1}} \dd \sigma,
\end{equation}
 if and only if $(M\smallsetminus\Omega, g)$ is isometric to
\begin{align}\label{eq:metric_cone}
    \left( [\rho_0, +\infty) \times \partial \Omega, \dd \rho \otimes \dd \rho + \left( \frac{\rho}{\rho_0}\right)^2 g^{}_{\partial \Omega}\right)&& \text{where }\rho_0 =\left(\frac{\abs{\partial \Omega}}{\AVR(g)\abs{\S^{n-1}}}\right)^{\frac{1}{n-1}}.
\end{align}
\end{theorem}

Some comments are in order about the above statements. First, we recall for the reader's convenience that the asymptotic volume ratio of $(M,g)$ is given by 
\begin{equation}
\AVR(g)= \lim_{r \to +\infty} \frac{\abs{B(o,r)}}{r^n \abs{\mathbb{B}^n}} \,,
\end{equation}
for some $o \in M$. The fact that, on complete manifolds with nonnegative Ricci curvature, the above limit is well defined and does not depend on the base point $o$, is a consequence of the classical Bishop-Gromov Volume Comparison Theorem. Moreover, one has that $0 \leq \AVR(g)\leq 1$, with $\AVR(g)=1$ if and only $(M,g)$ is the standard $n$-dimensional Euclidean space. Beside the intrinsic fundamental role played by manifolds with nonnegative Ricci curvature with Euclidean Volume Growth in geometric analysis, this class includes a diversity of explicit manifolds naturally arising from different fields, such as Asymptotically Locally Euclidean spaces (ALE for short) {\em gravitational instantons}. These are noncompact hyperkh\"aler Ricci Flat $4$-dimensional manifolds playing a  role in the study of Euclidean Quantum Gravity Theory, Gauge Theory and  String Theory (see~\cite{Hawking1977,Eguchi1979,Kronheimer1989,Kronheimer1989a,Minerbe2009,Minerbe2010,Minerbe2011}).

It is worth noticing that inequality~\eqref{eq:EMI_intro} is sharp, and in particular it provides the optimal Minkowski Inequality on manifolds with nonnegative Ricci curvature for \emph{outward minimising subsets}, see \cref{corollary-minimising}. These subsets are mean-convex and satisfy $\abs{\partial \Omega^*} = \abs{\partial \Omega}$, so that the Minkowski Inequality reads
\begin{equation}
\left(\frac{\abs{\partial \Omega}}{\abs{\S^{n-1}}}\right)^{\kern-.05cm\frac{n-2}{n-1}}\kern-.1cm\AVR(g)^{\frac{1}{n-1}}\leq \frac{1}{\abs{\S^{n-1}}}\int\limits_{\partial \Omega} { \frac{\HH}{n-1}} \dd \sigma,
\end{equation}
in this case. In addition to the Euclidean spaces, where it is immediately seen that balls achieve the identity in~\eqref{eq:EMI_intro}, the sharpness of this inequality is checked in a way greater generality. Namely, one can show that on manifolds that are just mildly asymptotically conical at infinity, one can find  a suitable sequence of bounded sets with smooth boundary saturating, in the limit, the above inequality.

Combining~\cref{thm:EMI_maintheorem} with the sharp Isoperimetric Inequality for manifolds with nonnegative Ricci curvature, firstly proved in dimension $3$ in~\cite[Theorem 1.4]{Agostiniani2018} and recently extended to any dimension by Brendle~\cite{Brendle2020} (see also \cite{Fogagnolo2020a,Johne2021,Balogh2021}), reading
\begin{equation}
\label{isoperimetric-intro}
   \frac{\abs{\Sf^{n-1}}^n}{\abs{\mathbb{B}^n}^{n-1}}\AVR(g)\leq\frac{\abs{\partial \Omega^*}^{n}}{\abs{\Omega^*}^{n-1}}
\end{equation}
we get the following sharp volumetric version of the Minkowski Inequality.
\begin{theorem}[Volumetric Minkowski inequality]
\label{volumetric}
Let $(M, g)$ be a complete Riemannian manifold with $\Ric \geq 0$ and Euclidean Volume Growth.
Then
\begin{equation}
\label{volumetricf}
\left(\frac{\abs{\Omega}}{\abs{\mathbb{B}^n}}\right)^{\!\!\frac{n-2}{n}} \!\! \AVR(g)^{\frac{2}{n}} \leq \frac{1}{\abs{\S^{n-1}}}\int\limits_{\partial \Omega} \abs{ \frac{\HH}{n-1}} \dd \sigma,
\end{equation}
where $\AVR(g)$ is the asymptotic volume ratio of $(M,g)$, $\HH$ is the mean curvature of $\partial \Omega$ with respect to the outward normal unit vector. Moreover, the equality is satisfied if and only if $(M,g)$ is isometric to the flat Euclidean space and $\Omega$ is a ball.
\end{theorem}
As for the Extended Minkowski Inequality, also~\eqref{volumetricf} is easily recognised to be sharp, while the rigidity statement directly follows from the rigidity of the Isoperimetric Inequality.
We finally point out that earlier contributions to the Volumetric Minkowski Inequality were given in~\cite{Chang2011} and~\cite{Qiu2015}, holding in the flat  Euclidean space and under stronger geometric assumptions on the boundary of $\Omega.$

\subsection{Outline of the proof.} We now describe the main features of our approach, that is in line with \cite{Agostiniani2016, Agostiniani2018, Fogagnolo2019, Agostiniani2019}. Given $(M,g)$ a Riemannian $n$-manifold, $n\geq 3$, with nonnegative Ricci curvature, and an open bounded subset $\Omega \subseteq M$ with smooth boundary we consider, 
for every $1<p<n$, the $p$-capacitary potential associated to $\Om$. This is the 
solution $u$ to the problem
\begin{equation}\label{pb-intro}
\begin{cases}[rcl@{\ \ }l]
\Delta^{\!(p)}_g u &=&0 &\text{ on $M \smallsetminus \overline{\Omega}$,}\\[.05cm]
u&=&1 & \text{ on $\partial\Omega$,}\\[.05cm]
u(x) &\to& 0 &\text{ as $\d(x , o)\to +\infty$,}
\end{cases}
\end{equation}
where $\Delta^{\!(p)}_g$ is the $p$-Laplace operator associated with the metric $g$, and $\d (\,\cdot\,, o)$ is the distance induced by $g$ to some fixed reference point $o$. Provided the manifold $(M,g)$ is $p$-nonparabolic (see \cref{def:pnop} below, as well as \cite{Holopainen1990,Holopainen1999}), the solution to problem~\eqref{pb-intro} exists and it is unique. Such a solution is commonly referred to as the \emph{$p$-capacitary potential} associated with $\Omega$. It is worth specifying that manifolds with Euclidean Volume Growth (i.e. $\AVR (g) > 0$) do satisfy the $p$-nonparabolicity assumption for $1 < p < n$ by   the characterisation given in \cite[Proposition 5.10]{Holopainen1999}. As a crucial step in our method, we will establish families of {\em monotonicity formulas}, holding along the level sets of the $p$-capacitary potentials associated with $\Om$. More precisely, for every $t \in [1,+\infty)$, we set 
\begin{align}\label{eq:Fbetaintro}
F^\beta_p(t)= t^{\beta\frac{(n-1)(p-1)}{(n-p)}}\kern-.5cm\int\limits_{\set{u=1/t}} \abs{\D u}^{(\beta+1)(p-1)}  \dd \sigma
\end{align}
and we show that for $\beta> (n-p)/[(p-1)(n-1)]$ the above quantity admits a nonincreasing $\CC^1(1,+\infty)$ representative.

Some remarks are mandatory at this stage. First of all, let us point out that the monotonicity statement provided here for the functions $F^\beta_p$'s holds in full generality and with no restriction on the geometry of $\Omega$. As such, it is new also for domains sitting in $\mathbb{R}^n$, where the same conclusions were provided in~\cite{Fogagnolo2019} only for convex domains, and in fact for smooth level sets flows. In the general case, it is well known that the level sets flow of $p$-harmonic functions might present a much less regular behaviour, since no general bound is available for the Hausdorff dimension of the critical set. To overcome these difficulties, the authors in~\cite{Agostiniani2019} settled for the \emph{effective inequalities}
\begin{align}\label{eq:effective_monotonicity}
\lim_{t\to + \infty} F^\beta_p(t)\leq F^\beta_p(1) && \text{and} && (F^\beta_p)'(1) \leq 0 \, .
\end{align}
The derivation of these two bounds, however, heavily relied on the compactness of the critical set of $u$, that is a particular feature of spaces with finite topology, and as such it is not directly viable in our setting (see \cite{Menguy2000}).   
In contrast with this, the present treatment provides the desirable extension to the nonlinear setting and to the general framework of nonnegatively Ricci curved $p$-nonparabolic manifolds of the monotonicity formulas discovered in~\cite{Colding2012, Colding2014a, Agostiniani2016, Agostiniani2018} for harmonic functions. As a second remark, to let the reader appreciate the $\CC^1$-regularity result, we observe that in principle even the fact that formula~\eqref{eq:Fbetaintro} yields a well posed definition is not granted for free. The most serious difficulty here is that the set of singular values cannot be controlled through the Sard's Theorem, since $p$-harmonic functions only enjoy a mild -- though optimal -- $\CC^{1,\beta}$-regularity.
We managed to solve these problems also taking advantage of recent insights given in \cite{Gigli2021}. The full statement of the Monotonicity Theorem is found in~\cref{prop:Monotonicity_Fbeta_sequence} below. 

Through the monotonicity of $F_p^\beta$, with $\beta = (p-1)^{-1}$, we arrive at the following \emph{$L^p$-Minkowski Inequality}
\begin{equation}
\label{lp-mink-intro}
\Cn_p(\Omega)^{\frac{n-p-1}{n-p}}\AVR(g)^{\frac{1}{n-p}}\leq \frac{1}{\abs{\S^{n-1}}} \int\limits_{\partial\Omega} \abs{\frac{\HH}{n-1}}^p \dd \sigma \, ,
\end{equation}
where $\Cn_p(\Omega)$ is the normalised $p$-capacity of $\Omega$ defined in~\eqref{eq:normalisedp-cap} below.
A major advantage we draw out of the full monotonicity of $F_p^\beta$ is the bypassing of the computation of its limit as $t \to + \infty$ when reaching for \eqref{lp-mink-intro}.  Indeed, this step is now replaced by a suitable contradiction argument
that combines the full monotonicity of our quantities with the sharp \emph{Iso-$p$-capacitary Inequality} (see \cref{thm:isopcapacitary} below)
\begin{equation}\label{eq:isopcap-intro}
     \frac{\Cn_p(\mathbb{B}^n)^n}{\abs{\mathbb{B}^n}^{n-p}}\AVR(g)^{p}\leq \frac{\Cn_p(\Omega)^n}{\abs{\Omega}^{n-p}}  .
\end{equation}
Such a statement is of independent interest in our opinion and can be achieved taking advantage of the already mentioned sharp Isoperimetric Inequality in manifolds with nonnegative Ricci curvature and Euclidean Volume Growth, following rather classical arguments (see e.g. \cite{Jauregui2012}).

With the $L^p$-Minkowski Inequality \eqref{lp-mink-intro} at hand, the Extended Minkowski Inequality \eqref{lp-mink-intro} simply follows by letting $p \to 1^+$ since 
\begin{equation}
\lim_{p \to 1^+} \Cn_p(\Omega)= \frac{\abs{\partial \Omega^*}}{\abs{\S^{n-1}}},
\end{equation}
as proven in \cite[Theorem 1.2]{Fogagnolo2020a}. This particular feature of our approach, namely the fact that the Minkowski Inequality is obtained as the limit of its $L^p$-versions, makes the rigidity statement a particularly nontrivial task, although we show that \eqref{lp-mink-intro} holds with an equality sign only on cones. This leads us to prove the rigidity statement, \cref{thm:rigidity}, through an  argument involving the study of the IMCF starting at boundaries of domains that saturate the Minkowski Inequality \eqref{eq:ExtendedMinkowskiIntro}. More precisely, we first show that for a short time the flow is smooth and given by constantly mean-curved totally umbilical hypersurfaces. This crucially exploits the nonnegativity of the Ricci curvature (\cref{lem:CMC_lemma}). Then, a splitting procedure along such flow, inspired by \cite{Huisken2001}, shows that an outer neighbourhood of $\partial \Omega$ is isometric to a truncated cone with the same volume ratio as $\AVR(g)$, and this allows to conclude (\cref{lem:BG_conico}).

\subsection{Further Monotonicity-Rigidity results}
Beside the Monotonicity-Rigidity properties of $F_p^\beta$ discussed above, we also establish analogous ones for the function
\begin{equation}
F^\infty_p(t)= t^{\frac{n-1}{n-p}}\sup_{\set{u =1/t}} \abs{\D u}.
\end{equation}

This is the content of \cref{prop:Monotonicity_Finfty_sequence}, that is again proved in the general setting of $p$-nonparabolic manifolds with nonnegative Ricci curvature, extending \cite[Theorem 1.3]{Fogagnolo2019}. As geometric consequences of this statement, we provide a rigidity result under pinching conditions and a sphere theorem for smooth boundaries in manifolds with $\ric \geq 0$ (see \cref{thm:p-pinching,thm:gradient pinching} below) and Euclidean Volume Growth. It is worth mentioning that the monotonicity of $F_p^\infty$ also leads to a new insight on the critical set of the $p$-capacitary potential, which we believe deserves some further investigations. Namely, it turns out that every level set of $u$ displays some nonempty relatively open region, where $\D u$ does not vanish, and where in particular $u$ is smooth (see  \cref{cor:regualirity_of_p-harmonic}).

\subsection{Summary}
In \cref{p-harmonic-sec} we report, for the ease of the reader, some relevant facts from the theory of $p$-harmonic functions on Riemannian manifolds, focusing on the regularity theory as well as on the existence and uniqueness of solutions to~\eqref{pb-intro}. Some important -- though already well known -- estimates and identities are also recalled in this section. \cref{sec2} is devoted to the proof of Monotonicity-Rigidity Theorems (see \cref{prop:Monotonicity_Fbeta_sequence,prop:Monotonicity_Finfty_sequence}). After having introduced a convenient conformally related setting, we restate them in this framework and we conclude the section with their proofs. In \cref{sec:consequences}, after having provided \eqref{eq:isopcap-intro}, we make use of these tools to prove the $L^p$-Minkowski Inequality (see \cref{thm:Lp_mink}), deduce the Extended Minkowski Inequality \cref{thm:EMI_maintheorem} and some rigidity results under pinching conditions as consequences of the Monotonicity-Rigidity Theorems.

\medskip

\noindent\textbf{Acknowledgements.} 
\emph{
The authors thank V. Agostiniani, L. De Masi, G. De Philippis, C. Mantegazza, A. Malchiodi, F. Oronzio and M. Vedovato for their interest in the present work as well as for precious discussions and comments during the preparation of this manuscript. The authors are members of Gruppo Nazionale per l'Analisi Matematica, la Probabilit\`a e le loro Applicazioni (GNAMPA), and they are partially funded by the GNAMPA project ``Aspetti geometrici in teoria del potenziale lineare e nonlineare".}

\section{The \texorpdfstring{$p$}{p}-capacitary potential in Riemannian manifolds}
\label{p-harmonic-sec}
We have collected here, for the sake of future reference, some substantially well known results that will be repeatedly applied in our arguments.
Before considering the specific case of problem \eqref{pb-intro}, we recall the definition of $p$-harmonic functions, as well as their regularity estimates. We then analyse the existence and uniqueness of solution $u_p$ to \eqref{pb-intro} on complete Riemannian manifolds. It turns out that these questions are intimately related to the notion of $p$-nonparabolicity. $p$-nonparabolic manifolds will then constitute the natural setting for the Monotonicity-Rigidity Theorems. We afterwards recall some global standard estimates on $u_p$ and its gradient as well as a Kato-type identity for $p$-harmonic functions.

\subsection{\texorpdfstring{$p$}{p}-harmonic functions and regularity}
\label{reg-subsec}
Given an open subset $U$ of a complete Riemannian manifold $(M, g)$, we say that $v \in W^{1, p}(U)$ is $p$-harmonic if 
\begin{equation}
\label{def-p-harm}
\int\limits_U \left\langle\abs{\D v}^{p-2} \D v \, \vert \, \D \psi\right\rangle \dd \mu =0.
\end{equation}
for any test function $\psi \in \CC_c^\infty(U)$. With $\langle\, \cdot \, \vert \, \cdot\, \rangle$ we denote as usual the scalar product induced by the underlying Riemannian metric $g$ on the tangent space at each point. Regularity results for $p$-harmonic functions (see \cite{Tolksdorf1984,DiBenedetto1983,Lieberman1988}) ensure that $v$ belongs to $\CC^{1, \beta}_{\loc}(U)$ for some $\beta\in (0,1)$ and is smooth around each point where $\abs{\D v}>0$.

Since the $\CC^{1,\beta}$-regularity is not sufficient to employ Sard's Theorem, we are going to heavily rely  on the coarea formula. We report it here for ease of further references. The statement below follows from \cite[Lemma 18.5 and Theorem 18.1]{Maggi2012} coupled with standard approximation results.

\crefalias{enumi}{proposition}
\begin{proposition}[Coarea formula]\label{prop:coarea}
Let $(M,g)$ be a complete Riemannian manifold. Consider a locally Lipschitz function $v:U \to [0,+\infty)$ on some open subset $U\subseteq M$, such that $v^{-1}([a,b])$ is compact for every $[a,b]\subset (0,+\infty)$. Then the following hold:
\begin{enumerate}[ref=\theproposition~(\arabic*),leftmargin=\parindent+6pt,align=left,labelwidth=\parindent,labelsep=6pt]
\item $\abs{\set{v=t}\cap \Crit (v)}=0$ for almost every $t \in [0,+\infty)$.
\item\label{item:coareaL1} for every measurable $f$ such that $f\abs{\D v}\in L^1_{\loc}(U)$ the function $f\in L^1(\set{v=t})$ for almost every $t\in(0,+\infty)$ and
\begin{equation}\label{eq:enached_coarea}
    \int\limits_{U} \psi(v) f \abs{ \D v} \dd \mu =\int\limits_I\psi(t)\int\limits_{\set{v=t}} f \dd \sigma \dd t,
\end{equation}
for every $\psi$ bounded measurable function compactly supported in $(0,+\infty)$. In particular, the function
\begin{equation}
    t \mapsto \int\limits_{\set{v=t}} f \dd \sigma \in L_{\loc}^1(0,+\infty)
\end{equation}
and its equivalence class does not depend on the representative of $f$.
\end{enumerate}
\end{proposition}

\begin{remark} \label{rmk:coarea_for_fraction}
If $h\in L^1_{\loc}(U)$ and $h=0$ almost everywhere on $\Crit(v)$, the function $f=h\abs{\D v}^{-1}$, satisfies the assumptions of \cref{item:coareaL1}. Clearly, if $f \in L^1(U)$,  \eqref{eq:enached_coarea} holds for every $\psi$ bounded measurable, even without compact support.
\end{remark}

With the idea of applying the previous result for $f=\abs{\D \abs{\D v}^{p-1}}$, a higher integrability degree of $p$-harmonic functions is required. We refer the reader to \cite[Lemma 2.1]{Lou2008} for a self-contained proof of the following lemma.

\begin{lemma}\label{rmk:vanishing_of_gradient}
Let $(M,g)$ be a complete Riemannian manifold and $U\subseteq M$ be a open subset. Given $v\in W^{1,p}(U)$ a $p$-harmonic function, then $\abs{\D v}^{p-1}\in W^{1,2}_\loc(U)$.
\end{lemma}

Given $U\subseteq M$ with Lipschitz boundary, a $p$-harmonic function $u\in W^{1,p}(U)$ attains some Dirichlet data $g\in L^p(\partial U)$ if $u$ coincides with $g$ on $\partial U$ in the sense of the trace operator.

%

\subsection{\texorpdfstring{$p$}{p}-nonparabolic manifolds and the \texorpdfstring{$p$}{p}-capacitary potential}
Given a noncompact Riemannian manifold $M$, we consider the $p$-capacitary potential of a bounded set with smooth boundary $\Omega \subset M$, that is a function $u \in W^{1, p}(M \smallsetminus \overline{\Omega})$ solving \eqref{pb-intro}. The function $u$ belongs to $\CC^{1,\beta}(M\smallsetminus \Omega)$ (see \cite{Lieberman1988}) and it is smooth near the points where the gradient does not vanish. In particular, by Hopf Maximum Principle in \cite[Proposition 3.2.1]{Tolksdorf1983} the datum on $\partial \Omega$ is attained smoothly.

We now focus on some classical sufficient conditions to ensure the existence of the $p$-capacitary potential, which turns out to be related to the notion of $p$-Green's function we are going to recall.

\begin{definition}[$p$-Green's function]
Let $(M, g)$ be a complete Riemannian manifold. Let $\mathrm{Diag}(M) = \{(x, x) \in M \times M \, \vert \, x \in M\}$. For $p \geq 1$, we say that $G_p : M \times M \smallsetminus \mathrm{Diag}(M) \to \R$ is a $p$-Green's function for $M$ if
it weakly satisfies $\Delta_p G(o, \,\cdot\,) = -\delta_o$ for any $o \in M$, where $\delta_o$ is the Dirac delta centred at $o$, that is, if it holds
\begin{equation}
\label{green-def}
\int\limits_M \Big\langle \abs{\D \, G_p{(o,\, \cdot\,)}}^{p-2} \, \D \,G_p{(o, \,\cdot\,)} \, \Big\vert \, \D \psi\Big\rangle \dd\mu =  \psi (o)
\end{equation}
for any $\psi \in \CC^{\infty}_c (M)$.
\end{definition}
The notion of $p$-Green's function calls for that of $p$-nonparabolic Riemannian manifold.

\begin{definition}[$p$-nonparabolicity]\label{def:pnop}
We say that a complete noncompact Riemannian manifold $(M, g)$ is $p$-nonparabolic if there exists a \emph{positive} $p$-Green's function $G_p : M \times M \smallsetminus \mathrm{Diag}(M) \to \R$. With the expression \emph{$p$-Green function} we are in fact referring to the positive minimal one.
\end{definition}
The notion of $p$-nonparabolicity is intimately related to existence of a solution to \eqref{pb-intro}, in that if the positive $p$-Green's function of a $p$-nonparabolic Riemannian manifold vanishes at infinity, then such solution exists for any open bounded subset $\Omega \subset M$ with smooth boundary. A complete and self contained proof of this fact is provided in the Appendix of \cite{Fogagnolo2020a}. We report the statement of such basic thought fundamental result.
\begin{theorem}[Existence of the $p$-capacitary potential]\label{thm:existence_of_p_potential}
Let $(M, g)$ be a complete noncompact $p$-nonparabolic Riemannian manifold. Let $\Omega \subset M$ be an open bounded subset with smooth boundary. Assume also that the $p$-Green's function $G_p$ satisfies $G_p(o, x) \to 0$ as $\d(o, x) \to +\infty$ for some $o \in M$. Then, there exists a unique solution $u_p$ to \eqref{pb-intro}. 
\end{theorem}

We find convenient to recall here the definition of $p$-capacity of an open bounded subset $\Omega \subset M$ together with a normalised version of it which turns out to be more advantageous for our computations.
\begin{definition}[$p$-capacity and normalised $p$-capacity]\label{def:capandnormpcap}
Let $(M, g)$ be a complete noncompact Riemannian manifold, and let $\Omega$ be an open bounded subset of $M$.
The $p$-capacity of $\Omega$ is defined as
\begin{equation}\label{cap}
\mathrm{Cap}_p (\Omega) = \inf \set{\,\int\limits_M \abs{\D v}^p \dd \mu\st v \in \CC_c^\infty(M),\, v \geq 1\text{ on } \Omega}.
\end{equation}
On the other hand, the normalised $p$-capacity of $\Omega$ is defined as
\begin{equation}\label{eq:normalisedp-cap}
\Cn_p(\Omega)=\frac{1}{\abs{\S^{n-1}}}\left(\frac{p-1}{n-p}\right)^{p-1}\capa_p(\Omega).
\end{equation}
\end{definition}

A function $u$ solving \eqref{pb-intro} realises the $p$-capacity of the initial set $\Omega$, and actually one can also characterise such quantity with a suitable integral on $\partial \Omega$. We resume these facts in the following statement.

\begin{proposition}
\label{cap-u-prop}
Let $(M, g)$ be a complete noncompact $p$-nonparabolic Riemannian manifold, for some $p > 1$. Let $\Omega \subset M$ be an open bounded subset with smooth boundary. Then the solution $u_p$ to \eqref{pb-intro} realises
\begin{equation}
\label{cap-u}
\Cn_p(\Omega) =\frac{1}{\abs{\S^{n-1}}} \left(\frac{p-1}{n-p}\right)^{p-1}  \int\limits_{M \smallsetminus \overline{\Omega}} \abs{\D u_p}^p \dd \mu
\end{equation}
Moreover, we have that
\begin{equation}
\label{p-cap-u}
\Cn_p(\Omega) =\frac{1}{\abs{\S^{n-1}}} \left(\frac{p-1}{n-p}\right)^{p-1}  \int\limits_{\set{u_p=1/t}} \abs{\D u_p}^{p-1} \dd\sigma.
\end{equation}
holds for almost every $t\in [1,+\infty)$, including any $1/t$ regular value for $u_p$. 
\end{proposition}

\proof The function $u_p$ can be approximated in $W^{1,p}(M\smallsetminus \overline{\Omega})$ by functions $\varphi$ in $\CC^{\infty}_c(M)$ which satisfies $\varphi\geq 1 $ on $\Omega$, then
\begin{equation}
    \capa_p(\Omega)\leq \int\limits_{M \smallsetminus \overline{\Omega}} \abs{ \D u_p}^{p} \dd \mu.
\end{equation}
On the other hand, the weak formulation in \eqref{def-p-harm} can be relaxed in duality with functions in $W^{1,p}_0(M\smallsetminus \overline{\Omega})$. Hence, taking any competitor $\psi\in\CC^\infty_c(M)$ with $\psi \geq 1$ on $\Omega$, $u_p- \psi \in W^{1,p}_0(M\smallsetminus \overline{\Omega})$ we get that 
\begin{align}
    \int\limits_{M \smallsetminus \overline{\Omega}} \abs{\D u_p}^p \dd \mu &=\int\limits_{M \smallsetminus \overline{\Omega}} \ip{ \abs{ \D u_p}^{p-2} \D u_p , \D u_p} \dd \mu =\int\limits_{M \smallsetminus \overline{\Omega}} \ip{ \abs{ \D u_p}^{p-2} \D u_p , \D \psi} \dd \mu.
\end{align}
Applying H\"{o}lder inequality to the right hand side, we are left with
\begin{equation}
     \int\limits_{M \smallsetminus \overline{\Omega}} \abs{\D v}^p \dd \mu\leq  \int\limits_{M \smallsetminus \overline{\Omega}} \abs{\D \psi}^p \dd \mu
\end{equation}
for every competitor $\psi$ in \eqref{cap}, proving \eqref{cap-u}. Since $\abs{\D u_p} \in L^p(M \smallsetminus \overline{\Omega})$, applying the coarea formula \eqref{eq:enached_coarea} with $f=\abs{\D u_p}^{p-1}$ to \eqref{cap-u} (see \cref{rmk:coarea_for_fraction}) one can obtain that
\begin{equation}\label{eq:capcoareaglobal}
    \capa_p(\Omega)=  \int\limits_0^1\int\limits_{\set{u_p=\tau}} \abs{\D u_p}^{p-1}\dd \sigma \dd \tau.
\end{equation}
Employing again the coarea formula \eqref{eq:enached_coarea} with $f=\abs{\D u_p}^{p-1}$ and the integration by part we get
\begin{align}
    \int\limits_0^1 \varphi'(\tau)\kern-.2cm \int\limits_{\set{u_p=\tau}} \abs{\D u_p}^{p-1} \dd \sigma \dd \tau  &=\int\limits_{M \smallsetminus \overline{\Omega}} \varphi'(u_p) \abs{ \D u_p}^p \dd \mu=-\int\limits_{M \smallsetminus \overline{\Omega}} \abs{ \D u_p}^{p-2} \ip{\D u_p, \D(\varphi(u_p))} \dd \mu \\& =\int\limits_{M \smallsetminus \overline{\Omega}} \varphi(u_p) \div( \abs{\D u_p}^{p-2} \D u_p) \dd \mu =0
\end{align}
for every $\varphi\in \CC^{\infty}_c (0,1)$, which gives that
\begin{equation}
    \tau \mapsto \int\limits_{\set{u_p=\tau}}  \abs{\D u_p}^{p-1} \dd \sigma 
\end{equation}
admits a constant representative, that coupled with \eqref{eq:capcoareaglobal} yields \eqref{p-cap-u}, denoting $t=1/\tau$. \endproof

In particular, evaluating \eqref{p-cap-u} at $t=1$, which is a regular value by the Hopf Maximum Principle \cite[Proposition 3.2.1]{Tolksdorf1983}, we have that
\begin{equation}
\label{cap-boundary}
\Cn_p (\Omega) = \frac{1}{\abs{\S^{n-1}}} \left(\frac{p-1}{n-p}\right)^{p-1}  \int\limits_{\partial \Omega} \abs{\D u_p}^{p-1} \dd\sigma.
\end{equation}

Moreover, one can actually relate the capacity of $\Omega_t= \set{u>1/t}\cup \Omega$ to the capacity of $\Omega$. The proof of the following lemma is contained in \cite[Lemma 3.8]{Holopainen1990}.

\begin{proposition}
Let $(M, g)$ be a complete noncompact $p$-nonparabolic Riemannian manifold, for some $p > 1$. Let $\Omega \subset M$ be an open bounded subset with smooth boundary. Then the solution $u_p$ to \eqref{pb-intro} realises
\begin{equation}\label{eq:scaling-invariant-cap}
    \Cn_p(\Omega_t)= t^{p-1}\Cn_p(\Omega)
\end{equation}
for every $t \in [1,+\infty)$, where $\Omega_t= \set{u>1/t}\cup \Omega$. In particular, the map $t\mapsto \Cn_p(\Omega_t)$ is smooth.
\end{proposition}

\subsection{Li-Yau-type estimates}
We provide a sharp lower estimate for the $p$-Green's function, extending the well known
\begin{equation}
\label{2-green-estimate}
\d(o, x)^{2 - n} \leq G_2(o, x)
\end{equation}
holding true for any couple of points $o, x$ belonging to a $2$-nonparabolic Riemannian manifolds with nonnegative Ricci curvature.
The proof of \eqref{2-green-estimate} builds on the Laplacian Comparison, that applies to show that
\[
\Delta \, \d(o,\, \cdot\,)^{2 - n} \geq 0
\] 
\emph{in the sense of distributions}. This amounts to say that
\begin{equation}
\label{lapl-distr}
- \int\limits_M \left\langle\D \, \, \d(o,\, \cdot\,)^{2-n} \, \big\vert \, \D\psi \right\rangle \dd\mu =  \int\limits_M  \d(o, \,\cdot\,)^{2-n} \Delta\psi \dd\mu \, \geq \, 0
\end{equation}
for any test function $\psi \in \mathscr{C}^{\infty}_c (M)$.
This leads to \eqref{2-green-estimate} substantially through the Maximum Principle. We refer the reader to \cite[Lemma 2.12]{Agostiniani2018} for details.
The nonlinear version of \eqref{2-green-estimate}, that, to our knowledge, has not been explicitly pointed out in literature yet, actually relies on \eqref{lapl-distr} too.
\begin{proposition}[Sharp lower bound for the $p$-Green's function]
\label{sharp-lower-prop}
Let $(M, g)$ be a complete $p$-nonparabolic Riemannian manifold with $\ric \geq 0$, for some $p > 1$. Let $o \in M$. Then, we have  
\begin{equation}
\label{sharp-lower-G}
\d (o, x)^{- \frac{n-p}{p-1}} \leq G_p(o, x)
\end{equation}
for any $x \in M \smallsetminus \{o\}$.
\end{proposition}
\begin{proof}
Fix for simplicity $o \in M$, and let $r(x) = \d(o, x)$.
We first show that $\Delta_p r^{- {(n-p)}/{(p-1)}} \geq 0$ holds in the weak sense, that is
\begin{equation}
\label{weak-subsol}
\int\limits_M \ip{ \abs{\D r^{- \frac{n-p}{p-1}}}^{p-2} \D r^{- \frac{n-p}{p-1}} \, , \, \D \psi } \dd\mu \leq 0
\end{equation}
for any $\psi \in \mathscr{C}^{\infty}_c(M)$. In fact, we have
\begin{equation}
\label{conto-subsol}
\begin{split}
\int\limits_M \ip{\abs{\D r^{- \frac{n-p}{p-1}}}^{p-2} \D r^{- \frac{n-p}{p-1}} , \D \psi } \dd\mu &= - \left(\frac{n-p}{p-1}\right)^{p-1} \int\limits_M r^{1-n}\ip{ \D  r ,\D \psi} \dd\mu \\
&= \frac{1}{n-2} \left(\frac{n-p}{p-1}\right)^{p-1}\int\limits_M \ip{ \D  r^{2-n} , \D \psi} \dd\mu \leq 0
\end{split}
\end{equation}
where the last inequality is the Laplacian Comparison Theorem \eqref{lapl-distr}.

Let now be $\delta > 0$. Since both $r^{- {(n-p)}/{(p-1)}}$ and $G_p$ vanish at infinity, we have  $r^{- {(n-p)}/{(p-1)}} \leq G_p + \delta$ on $\partial B(o, R)$ for any $R > 0$ big enough. On the other hand, the general result \cite[Theorem 12]{Serrin1964} ensures that $G_p(o, x)$ is asymptotic to $r(x)^{- {(n-p)}/{(p-1)}}$ as $\d(o, x) \to 0^+$, and thus we also get $r^{- {(n-p)}/{(p-1)}} \leq G_p + \delta$ on $\partial B(o, \epsilon)$ for any $\epsilon > 0$ small enough. Thus, applying the Comparison Principle to the subsolution  $r^{- {(n-p)}/{(p-1)}}$ and to the solution $G_p + \delta$ (with respect to the $p$-Laplacian), in the annulus $B(o, R) \smallsetminus \overline{ B(o, \epsilon)}$, we get
$r^{- {(n-p)}/{(p-1)}} \leq G_p + \delta$ on such annulus. Letting $\epsilon \to 0^+$ and $R \to +\infty$, we deduce that the same holds on the whole $M \smallsetminus \{o\}$. Finally, letting $\delta \to 0^+$, we are left with \eqref{sharp-lower-G}.
\end{proof}
We point out that \eqref{sharp-lower-G} is sharp, since the $p$-Green's function of $\R^n$ is exactly given by the formula $G_p(x, y) = \d (x, y)^{- {(n-p)}/{(p-1)}}$ for any $x \neq y$.
For what it concerns an upper bound for the $p$-Green's function, we observe that in \cite[Theorem 3.8]{Mari2019} it is shown, building on \cite[Proposition 5.10]{Holopainen1999}, that if $(M, g)$ in addition to the assumptions of \cref{sharp-lower-prop} has Euclidean Volume Growth, then we also have
\begin{equation}
\label{upper-G}
G_p(o, x) \leq \mathrm{C} \,  \d (o, x)^{- \frac{n-p}{p-1}} 
\end{equation}
for some constant $\mathrm{C}$ with well understood dependencies. 
In particular, coupling \eqref{sharp-lower-G} and \eqref{upper-G} with the Comparision Principle we deduce the following important estimate for the $p$-capacitary potential.
\begin{theorem}
\label{li-yau-u}
Let $(M, g)$ be a complete $p$-nonparabolic Riemannian manifold for some $p > 1$, with $\ric \geq 0$. Let $\Omega \subset M$ be a bounded subset with smooth boundary, and let $u_p$ be its $p$-capacitary potential. Then, there exists a positive constant $\Cn_1$ such that
\begin{equation}
\label{li-yauf}
\Cn_1\, \d (o, x)^{- \frac{n-p}{p-1}} \leq u_p(x) .
\end{equation}
for any $x \in M \smallsetminus \Omega$. 
If in addition $(M, g)$ has Euclidean Volume Growth, then there also exists another positive constant $\Cn_2$ such that  
\begin{equation}
\label{li-yauf2}
u_p(x) \leq \Cn_2\, \d (o, x)^{- \frac{n-p}{p-1}}.
\end{equation}
\end{theorem}
\begin{proof}
In light of \eqref{sharp-lower-G} and \eqref{upper-G}, this one holding true if $(M, g)$ satisfies the additional Euclidean Volume Growth assumption, it suffices to show that there exist positive constants $\mathrm{C}_1$ and $\mathrm{C}_2$ such that $\mathrm{C}_1 G_p \leq u_p \leq \mathrm{C}_2 G_p$. Choose any $\mathrm{C}_1 < 1 / \sup_{\partial \Omega} u_p$. Then, $\mathrm{C}_1 G_p < u_p$ on $\partial \Omega$. Moreover, since both $u_p$ and $G_p$ vanish at infinity, for any $\delta > 0$ we have $\mathrm{C}_1 G_p < u_p + \delta$ on $\partial B(o, R)$ for any $R$ big enough. The Comparison Principle applied to the $p$-harmonic functions $u_p + \delta$ and $G_p$ in $B(o, R) \smallsetminus \overline{\Omega}$ shows that $\mathrm{C}_1 G_p < u + \delta$ in the latter subset. The radius $R$ being arbitrarily big, this implies that, by passing to the limit as $R \to + \infty$, that $\mathrm{C}_1 G_p < u_p + \delta$  in the whole $M \smallsetminus \Omega$. Letting $\delta \to 0^+$ leaves with $\mathrm{C}_1 G_p \leq u_p$, and consequently with \eqref{li-yauf}. The inequality $u_p \leq \mathrm{C}_2 G_p$, yielding \eqref{li-yauf2}, is shown the same way.
\end{proof}
We now couple \eqref{li-yauf} with the general Cheng-Yau-type inequality for $p$-harmonic functions on manifolds with nonnegative Ricci curvature provided in \cite{Wang2010}. It asserts that a $p$-harmonic function $v$, with $p > 1$ defined in a ball $B(o, 2R) \subset M$, where $M$ is endowed with a Riemannian metric such that $\Ric \geq 0$, satisfies the estimate
\begin{equation}
\label{p-cheng-yau}
\sup_{B(o, R)}\frac{\abs{\D v}}{v} \leq \frac{\mathrm{C}}{R},
\end{equation} 
for a constant $\mathrm{C}$ depending only on the dimension of the ambient manifold and $p$.
With these tools we immediately obtain
\begin{proposition}
\label{gradient-boundu}
Let $(M, g)$ be a $p$-nonparabolic Riemannian manifold for some $p > 1$, with $\ric \geq 0$. Let $\Omega \subset M$ be a bounded subset with smooth boundary, and let $u_p$ be its $p$-capacitary potential. Then, there exists a positive constant $\mathrm{C}$ such that
\begin{equation}
\label{gradient-bounduf}
\abs{\D u_p}u_p^{-\frac{n-1}{n-p}} \leq \mathrm{C}
\end{equation}
holds on the whole $M \smallsetminus \Omega$.
\end{proposition}
\begin{proof}
By the $\mathscr{C}^{1}$-regularity of $u_p$, it clearly suffices to show that \eqref{gradient-bounduf} holds outside some compact set containing $\overline{\Omega}$. Let then $o \in \Omega$ and $R > 0$ be such that $\Omega \Subset B(o, R)$, and let $x \in M \smallsetminus \overline{B(o, 2R)}$. With this choice, we have $B(x, \d(o, x) - R) \Subset M \smallsetminus \overline{B(o, R)}$. Thus, applying inequality \eqref{p-cheng-yau} to the function $u_p$ in the ball $B(x, \d(o, x) - R)$ we get
\begin{equation}
\label{cheng-yau-applied1}
\frac{\abs{\D u_p}}{u_p^{\frac{n-1}{n-p}}} (x) \leq 2\mathrm{C} \frac{u_p(x)}{\d(o, x) - R} u_p^{-\frac{n-1}{n-p}}(x)\leq 4\Cn \frac{u_p^{-\frac{p-1}{n-p}}(x)}{\d(o, x)}
\end{equation}
and the rightmost hand side is bounded by means of \eqref{li-yauf}.
\end{proof}

\subsection{Kato-type identity and a warped product splitting theorem}
Finally we give the statement of the refined Kato-type identity for $p$-harmonic functions obtained in \cite[Proposition 4.4]{Fogagnolo2019}, that will be at the core of the monotonicity and rigidity of $F_p^\beta$. 

\begin{definition}[Geometry of level sets and orthogonal decomposition]\label{def:Geometry_and_orthogonal_decomposition}
Let $(M,g)$ be a Riemannian manifold and $v$ be a smooth function on $M$. At any point where $\abs{\D v}\neq 0$ we denote by $\hh$ and $\HH$ respectively the second fundamental form and the mean curvature of the level set of $u$ with respect to the unit normal $\D v / \abs{\D v}$ and $g^\top$ the metric induced by $g$ on the level set of $u$. Finally, for a given differentiable function $f$, we denote by $\D^\top f$ the tangential part of the gradient, according to the orthogonal decomposition
\begin{align}
\D^\perp f = \ip{\D f, \frac{ \D v}{\abs{ \D v}}}\frac{\D v}{\abs{ \D v}}&& \text{ and} && \D^\top f= \D f- \D^\perp f.
\end{align}
\end{definition}
In particular, the following formula holds
\begin{equation}\label{eq:orthogonal_decomposition}
\abs{ \D \abs{ \D f}}^2 = \abs{ \D^\top\abs{\D f}}^2 + \abs{ \D^\perp\abs{\D f}}^2.
\end{equation}
We are now ready to state the Kato-type identity for $p$-harmonic function.

\begin{proposition}[Kato-type identity] \label{prop:Kato-type}
Let $(M,g)$ be a Riemannian manifold and let $v$ be a $p$-harmonic function on some subset of $M$, $p>1$. Then, in an open neighbourhood of a point where $\abs{\D v}\neq 0$, the following identity holds
\begin{equation}\label{eq:Kato-type_identity}
\abs{ \D \D v}^2 - \left(1+ \frac{(p-1)^2}{n-1}\right) \abs{ \D \abs{\D v}}^2 =\abs{ \D v}^2 \abs{ \hh- \frac{\HH}{n-1}g^{\top}}^2 + \left( 1- \frac{(p-1)^2}{n-1}\right) \abs{ \D^\top \abs{ \D v}}^2,
\end{equation}
according to the notation in \cref{def:Geometry_and_orthogonal_decomposition}. Moreover, if for some $t_0 \in \R$, $\abs{\D v}>0$ and 
\begin{align}
\abs{ \hh- \frac{\HH}{n-1}g^{\top}}^2=0, &&\abs{ \D^\top \abs{ \D v}}^2=0
\end{align}
hold at each point of $\set{v \geq t_0}$, then the Riemannian manifold $(\set{v\geq t_0},g)$ is isometric to the warped product $([t_0,+\infty) \times \set{v =t_0}, \dd t \otimes \dd t + \eta^2(t) g^{}_{\set{v =t_0}})$, where the relation between $v$, $\eta$ and $t$ is given by
\begin{equation}\label{eq:Kato-type_rigidity}
\eta(t)= \left( \frac{v'(t_0)}{v'(t)}\right)^{\frac{p-1}{n-1}}.
\end{equation}
\end{proposition}

\section{Monotonicity-Rigidity Theorems}
\label{sec2}
In this section we are going to prove our \textit{monotonicity formulas} in the $p$-nonparabolic setting. The results we present here are the natural extensions of the ones shown in \cite{Agostiniani2016,Agostiniani2018} as well as of the ones obtained in
\cite{Fogagnolo2019,Agostiniani2019}. In the first two mentioned papers
the authors established the monotonicity in the case of the harmonic potential, respectively in $\R^n$ and in a general $2$-nonparabolic manifold with nonnegative Ricci curvature, whereas in the second two papers an analogous theory has been developed in the case of the $p$-capacitary potential in the Euclidean setting. More precisely, in \cite{Fogagnolo2019}, the authors worked out the smooth computations, and took advantage of the fact that the $p$-capacitary potential associated with a convex domain is smooth and has no critical points (see \cite{Colesanti2015,Lewis1977}), whereas the main technical achievement in \cite{Agostiniani2019} is the treatment of the general case, when the critical points are present and even possibly arranged in sets of full measure. On the other hand, the approach presented in \cite{Agostiniani2019} only produces \emph{effective inequalities} \eqref{eq:effective_monotonicity}, that are anyway sufficient to prove \cref{thm:EMI_maintheorem} in the flat setting, as mentioned in the Introduction. Here, we extend these results to the the setting of $p$-nonparabolic manifolds and we improve them, establishing the full monotonicity of the integral quantities defined in \eqref{eq:monotonefunction} along the $p$-capacitary level sets flow. 

As usual, the main difficulty amounts to ensure that the monotonicity survives the singular values of $u$, that, as far as we known, could even form a set of positive measure. Inspired by the analysis in \cite{Gigli2021}, where the authors were forced to face severe technical problems caused by the typical low regularity of the nonsmooth setting, we compute the derivative of our integral quantities \eqref{eq:monotonefunction} in the distributional sense, appealing to the full strength of the coarea formula in \cref{prop:coarea}, and exploiting the integrability properties of the $p$-harmonic functions in \cref{rmk:vanishing_of_gradient}.

\medskip

\noindent {\em From now on, unless where it is necessary, we fix $1<p<n$ and we drop the subscript $p$ when we consider a solution $u_p$ to the problem \eqref{pb-intro}.}

\medskip

\subsection{Statement of the Monotonicity-Rigidity Theorems}

Let $u:M\smallsetminus \Omega\to \R$ be a solution of \eqref{pb-intro}. For $\beta \in [0,+\infty)$ we consider the function
\begin{equation}\label{eq:monotonefunction}
F^\beta_p(t)= t^{\beta\frac{(n-1)(p-1)}{(n-p)}}\kern-.5cm\int\limits_{\set{u=1/t}} \abs{\D u}^{(\beta+1)(p-1)}  \dd \sigma
\end{equation}
defined for every $t\geq 1$ such that $\abs{\set{u=1/t}\cap \Crit (u)}=0$, which is fulfilled for almost every $t \in [1,+\infty)$ by \cref{prop:coarea}. We also set
\begin{equation}\label{eq:monotone_function_infty}
F^\infty_p(t)= t^{\frac{n-1}{n-p}}\sup_{\set{u =1/t}} \abs{\D u},
\end{equation}
that is defined on the whole $[1,+\infty)$. If $1/t$ is a regular value for $u$, then $F^\beta_p$ is differentiable at $t$ for every $\beta\in [0,+\infty)$ and its derivative is
\begin{align} \label{eq:derivative_of_F_standard}
(F^\beta_p)'(t)=-\beta t^{\beta \frac{(n-1)(p-1)}{(n-p)}-2}\kern-.4cm\displaystyle\int\limits_{\set{u=1/t}}\kern-.3cm\abs{\D u}^{(\beta+1)(p-1)-1}\kern-.1cm\left( \HH -{\textstyle\frac{(n-1)(p-1)}{(n-p)}} \abs{\D \log u}\right) \dd \sigma.
\raisetag{.5em}
\end{align}

As said before, the aim of this section is to prove Monotonicity-Rigidity Theorems for $t \mapsto F_p^\beta(t)$ and $t \mapsto F_p^\infty(t)$.

\begin{theorem}[Monotonicity-Rigidity theorem for $F_p^\beta$]\label{prop:Monotonicity_Fbeta_sequence}
Let $(M,g)$ be a $p$-nonparabolic Riemannian manifold with $\Ric \geq 0$. Let $\Omega\subseteq M$ be a bounded open subset with smooth boundary. Let $F_p^\beta$ be the function defined in \eqref{eq:monotonefunction} with $(n-p)/[(n-1)(p-1)]< \beta<+\infty$. Then $F^\beta_p$ belongs to $W^{2,1}(1,+\infty)$ and the identities
\begin{allowdisplaybreaks}
\begin{align}\label{eq:derivative_of_Fbeta}
\begin{split}(F_p^\beta)'(t)&= -\beta\left(\frac{(n-2)(p-1)}{(n-p)}\right)^{(\beta+1)(p-1)}\kern-1
cm \int\limits_{\set{u \leq 1/t}\smallsetminus \Crit(u)}\kern-.8cm  u^{2- \beta \frac{(p-1)(n-1)}{(n-p)}}\abs{\D u}^{(\beta+1)(p-1)-1}\\
&\kern.5cm\vphantom{\int\limits_{\set{u \leq 1/t}}}\left\lbrace\left[\beta - \frac{(n-p)}{(n-1)(p-1)}\right]  \left[\HH - \left[\frac{(n-1)(p-1)}{(n-p)}\right] \abs{\D \log u}\right]^2+ \abs{\hh - \frac{\HH}{n-1}g^\top}^2\right.\\
&\kern+3cm+(p-1)\left[\beta +\frac{p-2}{p-1}\right]\frac{\abs{\D^\top \abs{\D u}}^2}{\abs{\D u}^2}\vphantom{\abs{\hh - \frac{\HH}{n-1}g^\top}^2}\vphantom{\int\limits_{\set{u \leq 1/t}}}+\left.\Ric\left(\frac{\D u}{\abs{\D u}}, \frac{\D u}{\abs{\D u}}\right)\vphantom{ \left[\HH - \left[\frac{(n-1)(p-1)}{(n-p)}\right] \abs{\D \log u}\right]^2}\right\rbrace \dd \mu\vphantom{\int\limits_{\set{u \leq 1/t}}} \end{split}\kern.5cm
\intertext{and}
\label{eq:2nd_derivative_of_Fbeta}%
\begin{split}(F_p^\beta)''(t)&= \beta\, \left(\frac{(n-2)(p-1)}{(n-p)}\right)^{(\beta+1)(p-1)}t^{\beta\frac{(n-1)(p-1)}{(n-p)}-2}\kern-.2cm  \int\limits_{\set{u = 1/t}}\kern-.2cm  \abs{\D u}^{(\beta+1)(p-1)-2}\\&\kern.5cm\vphantom{\int\limits_{\set{u \leq 1/t}}}\left\lbrace\left[\beta - \frac{(n-p)}{(n-1)(p-1)}\right]  \left[\HH - \left[\frac{(n-1)(p-1)}{(n-p)}\right] \abs{\D \log u}\right]^2+ \abs{\hh - \frac{\HH}{n-1}g^\top}^2\right.\\
&\kern+3cm+(p-1)\left[\beta +\frac{p-2}{p-1}\right]\frac{\abs{\D^\top \abs{\D u}}^2}{\abs{\D u}^2}\vphantom{\abs{\hh - \frac{\HH}{n-1}g^\top}^2}\vphantom{\int\limits_{\set{u \leq 1/t}}}+\left.\Ric\left(\frac{\D u}{\abs{\D u}}, \frac{\D u}{\abs{\D u}}\right)\vphantom{ \left[\HH - \left[\frac{(n-1)(p-1)}{(n-p)}\right] \abs{\D \log u}\right]^2}\right\rbrace \dd \mu\vphantom{\int\limits_{\set{u \leq 1/t}}},\end{split}\kern.5cm
\end{align}
\end{allowdisplaybreaks}%
hold for almost every $t\in[1,+\infty)$. In particular, $F^\beta_p$ admits a convex and monotone nonincreasing $\CC^1$ representative. Moreover, $(F^\beta_p)'(t_0)=0$ at some $t_0\geq 1$ such that $1/t_0$ regular value for $u$ if and only if $(\set{u\leq 1/t_0}, g)$ is isometric to
\begin{align}
\left([\tau_0,+\infty) \times \set{u=1/t_0}, \dd \tau \otimes \dd \tau  + \left(\frac{\tau }{\tau_0}\right)^2 g^{}_{ \set{u=1/t_0}}\right),&& \text{where } \tau _0 = \left(\frac{\abs{\set{u=1/t_0}}}{\AVR(g)\abs{\S^{n-1}}}\right)^{\frac{1}{n-1}}.
\end{align}
In this case $\set{u=1/t_0}$ is a connected totally umbilical hypersurface with constant mean curvature in $(M\smallsetminus \Omega, g)$ .
\end{theorem}

We highlight that the rigidity statement is expressed in terms of the derivative. However, if $F_p^\beta(t)=F_p^\beta(T)$ for $1\leq t<T < +\infty$ such that $1/t$ and $1/T$ are regular value for $u$, the rigidity statement still triggers. Indeed, since the set of regular values is open, monotonicity ensures the existence of a decreasing sequence $(t_j)_{j\in \N}$ such that $t_j\to t$ as $j \to +\infty$, $1/t_j$ is regular for $u$ and $(F_p^\beta)'(t_j)=0$. Since $t\mapsto F_p^\beta(t)$ is smooth in a neighbourhood of $t$, this implies that $(F_p^\beta)'(t)=0$ hence the splitting of $\set{u\leq 1/t}$.

\begin{theorem}[Monotonicity-Rigidity theorem for $F^\infty_p$]\label{prop:Monotonicity_Finfty_sequence}
Let $(M,g)$ be a $p$-nonparabolic Riemannian manifold with $\Ric \geq 0$. Let $\Omega\subseteq M$ be a bounded open subset with smooth boundary. Let $F_p^\infty$ be the function defined in \eqref{eq:monotone_function_infty}. Then $F^\infty_p$ is a continuous monotone nonincreasing function. Furthermore, we have
\begin{equation}\label{eq:mean_curvature_inequality_background}
\left[\HH_g - \frac{(n-1)(p-1)}{(n-p)}\abs{\D \log u}_g\right](x_t)=-(p-1)\frac{\partial}{\partial \nu_t} \log \frac{ \abs{ \D u}_g}{u^{\frac{n-1}{n-p}}}(x_t)\geq 0
\end{equation}
where $x_t\in \set{u=1/t}$ is the point where $\sup_{\set{u=1/t}}\abs{\D u}_g/u^{(n-1)/(n-p)}$ is achieved and $\nu_t=-\D u / \abs{\D u}_g$ is the unit normal to $\set{u=1/t}$.
Moreover, $F_p^\infty(t_0)=F_p^\infty(T)$ for some $t_0<T$ or the equality holds in \eqref{eq:mean_curvature_inequality_background} for some $t_0$ such that $1/T$ and $1/t_0$ are regular for $u$ if and only if $(\set{u\leq 1/t_0}, g)$ is isometric to 
\begin{align}
\left([\tau_0,+\infty) \times \set{u=1/t_0}, \dd \tau \otimes \dd \tau  + \left(\frac{\tau }{\tau_0}\right)^2 g^{}_{ \set{u=1/t_0}}\right),&& \text{where } \tau _0 = \left(\frac{\abs{\set{u=1/t_0}}}{\AVR(g)\abs{\S^{n-1}}}\right)^{\frac{1}{n-1}}.
\end{align}
In this case $\set{u=1/t_0}$ is a connected totally umbilical hypersurface with constant mean curvature in $(M\smallsetminus \Omega, g)$.
\end{theorem}

A direct consequence of the monotonicity of $F^\infty_p$ is the following regularity theorem for the $p$-capacitary potential.

\begin{corollary}\label{cor:regualirity_of_p-harmonic}
The function $F^\infty_p$ is strictly positive. In particular, every level of $u$ has at least one regular point.
\end{corollary}

We want also to emphasise that these theorems can be applied in particular in $\R^n$ for every $\Omega$ open bounded with smooth boundary, where they naturally extend the Monotonicity-Rigidity Theorems in \cite{Fogagnolo2019, Agostiniani2019}.

We conclude this introduction rewriting the function $F_p^\beta$ and $F_p^\infty$ defined in \eqref{eq:monotonefunction} and \eqref{eq:monotone_function_infty} in a different formulation. We make use of this tool only to simplify computations, but as shown in \cite{Agostiniani2016,Fogagnolo2019,Agostiniani2020} Monotonicity-Rigidity theorems have their counterpart in this framework. Let $(M,g)$ be a complete $p$-nonparabolic Riemannian manifold with $\Ric\geq 0$ and $u:M\smallsetminus\Omega \to \R$ be the solution of the problem \eqref{pb-intro}. We consider the conformally related Riemannian manifold $(M\smallsetminus \Omega, \tilde{g})$, where $\tilde{g}$ is given by
\begin{equation}\label{eq:conformal_metric}
\tilde{g}=u^{2\left(\frac{p-1}{n-p}\right)} g.
\end{equation}
It is also convenient to consider the new variable
\begin{equation}\label{eq:definitionphi}
\varphi=-\frac{(p-1)(n-2)}{(n-p)} \log u,
\end{equation}
so that the metric $\tilde{g}$ can be equivalently rewritten as
\begin{equation}
\tilde{g}=\ee^{-\frac{2 \varphi}{n-2}}g.
\end{equation}

With the same formal computation as in \cite{Fogagnolo2019}, one can prove that $\Delta_{\tilde{g}}^p \varphi =0$ on $M\smallsetminus \overline{\Omega}$ where $\Delta_{\tilde{g}}^p$ is the $p$-Laplace operator with respect to the metric $\tilde{g}$.


\medskip

\noindent{\em From now on, given $(M,g)$ a $p$-nonparabolic manifold with $\Ric\geq0$ and $u$ a solution to \eqref{pb-intro}, $\varphi$ will be the function otteined by $u$ trough \eqref{eq:definitionphi} whereas $\tilde{g}$ will indicate the metric on $M\smallsetminus { \Omega}$ obtained from $u$ and $g$ trough \eqref{eq:conformal_metric}.}

\medskip

The gradient of $\phi$ is related to the one of $u$ by
\begin{equation}
\label{gradphi-relation}
\abs{\nabla \phi}_{\tilde{g}} = \frac{(n-2)(p-1)}{(n-p)} \frac{\abs{\D u}_g}{u^{\frac{n-1}{n-p}}},
\end{equation}
where $\nabla$ is the Levi-Civita connection associated to the metric $\tilde{g}$. We can observe that if $t$ is a regular value for $u$ then $ s= -[(p-1)(n-2)/(n-p)] \log t$ is a regular value for $\phi$, thanks to \eqref{eq:definitionphi} and the previous relation. Moreover, we recognise from the above expression and the estimate \eqref{gradient-bounduf} the fundamental property of $\abs{\nabla \phi}_{\tilde{g}}$ to be uniformly bounded, that is, there exists a constant $\mathrm{C}$ such that
\begin{equation}
\label{gradphi-bound}
\abs{\nabla \phi}_{\tilde{g}} \leq \mathrm{C}
\end{equation}
on the whole $M \smallsetminus \Omega$.

Using \eqref{gradphi-relation} the family of functions $t\mapsto F^\beta_p(t)$ for $\beta \in [0,+\infty]$ defined in \eqref{eq:monotonefunction} and \eqref{eq:monotone_function_infty} can be rewritten in terms of $\tilde{g}$ and $\varphi$ obtained through \eqref{eq:definitionphi} and \eqref{eq:conformal_metric}. For any $\beta \in [0,+\infty)$ we can now consider the function
\begin{equation}\label{eq:Phidef}
\Phi^\beta_p(s)=\int\limits_{\set{\varphi=s}} \abs{\nabla \varphi}^{(\beta+1)(p-1)}_{\tilde{g}}\dd \sigma_{\tilde{g}},
\end{equation}
whenever $s\geq 0$ is such that $\abs{\set{\varphi=s} \cap \Crit(\varphi)}$. Correspondingly we set 
\begin{equation}\label{eq:Phiinfty}
\Phi^\infty_p(s)= \sup_{\set{\varphi=s}}\abs{\nabla \phi}_{\tilde{g}}
\end{equation}
which is defined on the whole $[0,+\infty)$. The function $\Phi^\beta_p$ can be obtained from $F_p^\beta$ through a change of variable that is
\begin{equation}\label{eq:Phibeta_Fbeta}
\Phi^\beta_p(s)= F^\beta_p\left(\ee^{\frac{(n-p)}{(p-1)(n-2)}s}\right).
\end{equation}
For $\beta < +\infty$ it thus holds that
\begin{equation}\label{eq:dPhibeta_dFbeta}
(\Phi^\beta_p)'(s)= \frac{(n-p)}{(p-1)(n-2)}\ee^{\frac{(n-p)}{(p-1)(n-2)}s}(F^\beta_p)'\left(\ee^{\frac{(n-p)}{(p-1)(n-2)}s}\right)
\end{equation}
for almost every $s \in [0,+\infty)$. The previous relations reveal how proving the Monotonicity results for $F_p^\beta$ and $F_p^\infty$, stated in \cref{prop:Monotonicity_Fbeta_sequence} and \cref{prop:Monotonicity_Finfty_sequence}, are equivalent to show the same one for $\Phi^\beta_p$ and $\Phi^\infty_p$. The same argument applies for the regularity of $F^\beta_p$.

\subsection{Proof of Monotonicity-Rigidity Theorems} 
A basic property we will need is the essential uniform boundedness of $\Phi_p^\beta$ of $\Phi_p^{\infty}$ defined in \eqref{eq:Phidef} and \eqref{eq:Phiinfty}.

\begin{lemma}
\label{bound-Phi}
Let be $1 < p < n$, and $(M, g)$ be a $p$-nonparabolic Riemannian manifold. Let $\Omega \subset M$ be a open bounded subset with smooth boundary. For every $\beta\in[0,+\infty)$, $\Phi_p^\beta$ is essentially uniformly bounded, namely $\Phi_p^\beta(s) \leq C$ for almost every $s \in [0, +\infty)$, including any $s$ that is regular for $\phi$. Moreover, the function $\Phi_p^{\infty}$ is uniformly bounded. 
\end{lemma}

\proof
It suffices to write $\Phi^\beta_p$ as
\begin{align}
\Phi^\beta_p(s)&= \int\limits_{\set{\varphi =s}} \abs{ \nabla \varphi}_{\tilde{g}}^{(\beta+1)(p-1)} \dd \sigma_{\tilde{g} }
\leq \Cn^{\beta(p-1)}\int\limits_{\set{\varphi =s}} \abs{ \nabla \varphi}_{\tilde{g}}^{p-1} \dd \sigma_{\tilde{g}} \\&=\Cn^{\beta(p-1)}\left[\frac{(n-2)(p-1)}{(n-p)}\right]^{p-1}\int\limits_{\set{u =1/t}}\abs{\D u}^{p-1} \dd \sigma 
\end{align}
where $\mathrm{C}$ is the constant in \eqref{gradphi-bound}, the last identity is due to \eqref{gradphi-relation} and \eqref{eq:definitionphi} taking $ s= -[(p-1)(n-2)/(n-p)] \log t$. By \eqref{p-cap-u} we have that the integral on the rightmost hand side coincides with $\mathrm{Cap}_p (\Omega)$ for almost any $t$, including any of those such that $1/t$ is a regular value for $u$. This settles the boundedness of $\Phi_p^\beta$ for finite $\beta$. On the other hand the uniform boundedness of $\Phi^\infty_p$ is a direct consequence of \eqref{gradphi-bound} alone.
\endproof
From now on, we will drop the subscript $\tilde{g}$ whenever it is clear to which metric we are referring.

\smallskip

Suppose by now that $\beta\in [0,+\infty)$ and consider the vector field
\begin{equation}\label{eq:Xdef}
X=\ee^{-\frac{(n-p)}{(n-2)(p-1)}\varphi}\abs{\nabla \varphi}^{p-2} \left( \nabla\abs{\nabla \varphi}^{\beta(p-1)} + (p-2)\nabla^\perp \abs{\nabla \varphi}^{\beta(p-1)} \right),
\end{equation}
defined in a neighbourhood of each point such that $\abs{\nabla\varphi}\neq 0$. This vector field is related to derivative of $\Phi^\beta_p$, due to the following identity.

\begin{proposition}\label{prop:derivativeidentity}
Let $(M,g)$ be a $p$-nonparabolic Riemannian manifold with $\Ric\geq 0$. For every $\beta \in[0,+\infty)$, the function $s \mapsto\Phi^\beta_p(s)$ defined in \eqref{eq:Phidef} belongs to $W^{1,1}_{\loc}(0,+\infty)$ and its derivative is given by
\begin{equation}\label{eq:derivativeidentity}
\ee^{-\frac{(n-p)}{(n-2)(p-1)}s}(\Phi^\beta_p)'(s)= \frac{1}{p-1}\int\limits_{\set{\varphi=s}}\ip{X, \frac{\nabla \varphi}{\abs{\nabla \varphi}}} \dd \sigma.
\end{equation}
for almost every $s \in [0,+\infty)$, where $X$ is the vector field defined in \eqref{eq:Xdef}. 
\end{proposition}

\proof By the definition of $X$, it is easy to check that
\begin{equation}
\ee^{-\frac{(n-p)}{(n-2)(p-1)}\varphi}\ip{\abs{\nabla \varphi}^{p-2} \nabla \abs{\nabla \varphi}^{\beta(p-1)}, \frac{\nabla \varphi}{\abs{\nabla \varphi}}} = \frac{1}{p-1} \ip{X, \frac{\nabla \varphi}{\abs{\nabla \varphi}}},
\end{equation}
holds around each point such that $\abs{\nabla\varphi}\neq 0$. Hence, it remains only to prove that $\Phi_p^\beta (s) \in W^{1, 1}_{\mathrm{loc}} (0 + \infty)$ and that
\begin{equation}
(\Phi^\beta_p)'(s)= \int\limits_{\set{\varphi=s}} \ip{ \abs{\nabla \varphi}^{p-2} \nabla \abs{\nabla \varphi}^{\beta(p-1)} , \frac{\nabla \varphi}{\abs{\nabla \varphi}}} \dd \sigma
\end{equation}
holds for almost any $s \in (0, \infty)$.
Let $\eta \in \CC^\infty_c(0,+\infty)$. Since $\abs{\nabla \varphi}\leq C$ by \eqref{gradphi-bound}, applying the coarea formula \eqref{eq:enached_coarea} with $f= \abs{\nabla \varphi}^{(\beta+1)(p-1)}$ and the chain rule we obtain that
\begin{align}
\int\limits_0^{+\infty} \eta'(s)\Phi_p^\beta(s)\dd s &= \int\limits_0^{+\infty}\eta'(s)\kern-.2cm\int\limits_{\set{\varphi=s}}\kern-.2cm\abs{\nabla \varphi}^{(\beta+1)(p-1)}\dd \sigma\dd s=\int\limits_{M\smallsetminus\overline{\Omega}}\kern-.2cm\eta'(s)\ip{\nabla \varphi, \nabla \varphi}\abs{\nabla \varphi}^{(\beta+1)(p-1)-1}\dd \mu\\
&= \int\limits_{M\smallsetminus\overline{\Omega}}\ip{\nabla(\eta(\varphi)), \nabla \varphi}\abs{\nabla \varphi}^{(\beta+1)(p-1)-1}\dd \mu.
\end{align}
Integrating by parts the right hand side, $\Delta^{(p)} \varphi =0$ yield
\begin{equation}
\int\limits_0^{+\infty} \eta'(s)\Phi_p^\beta(s)\dd s  = -\int\limits_{M\smallsetminus \overline{\Omega}} \eta(\varphi) \ip{ \abs{\nabla \varphi}^{p-2} \nabla \abs{\nabla \varphi}^{\beta(p-1)} , \nabla \varphi}\dd \mu.
\end{equation}
Thanks to \eqref{gradphi-bound} and \cref{rmk:vanishing_of_gradient}, we are in position to apply the coarea formula in \cref{prop:coarea} with $f= \ip{\abs{\nabla \varphi}^{p-2} \nabla \abs{\nabla \varphi}^{\beta(p-1)} , \nabla \varphi}/\abs{\nabla \varphi}$ (see \cref{rmk:coarea_for_fraction}), to get
\begin{equation}
\int\limits_0^{+\infty} \eta'(s)\Phi_p^\beta(s)\dd s  = -\int\limits_0^1\eta(s) \int\limits_{\set{\varphi=s}}\ip{ \abs{\nabla \varphi}^{p-2} \nabla \abs{\nabla \varphi}^{\beta(p-1)} , \frac{\nabla \varphi}{\abs{\nabla \varphi}}} \dd \sigma\dd s,
\end{equation}
which ensures both that $\Phi_p^\beta\in W^{1,1}_{\loc}(0,+\infty)$ and \eqref{eq:derivativeidentity}.
\endproof


The nonnegative divergence of $X$ is what substantially rules the monotonicity of $\Phi_p^\beta$, and this is true when $\beta$ ranges in a suitable set of parameters.
\begin{lemma}[Divergence of $X$] \label{lem:divergence_of_X}
Let $(M,g)$ be a $p$-nonparabolic manifold and $X$ be the vector field defined in \eqref{eq:Xdef}. Then
\begin{equation}
\div X=\ee^{- \frac{(n-p)}{(n-2)(p-1)}\varphi} Q,
\end{equation}
holds at any point such that $\abs{\nabla \varphi}$, with
\begin{align}
\begin{split}
Q&=\beta(p-1)\abs{\nabla \varphi}^{\beta(p-1)+p-2}\left\lbrace\abs{\hh - \frac{\HH}{n-1}\tilde{g}^\top}^2 \right.+(p-1)\left[\beta+ \frac{p-2}{p-1}\right]\frac{\abs{\nabla^\top\abs{\nabla \varphi}}^2}{ \abs{\nabla \varphi}^2}\\\label{eq:Qdef}
&\kern+1cm+\left.(p-1)^2 \left[\beta - \frac{(n-p)}{(p-1)(n-1)}\right]\frac{\abs{\nabla^\perp\abs{\nabla \varphi}}^2}{ \abs{\nabla \varphi}^2}+\Ric_g\left(\frac{\nabla \varphi}{ \abs{\nabla \varphi}^2}, \frac{\nabla \varphi}{ \abs{\nabla \varphi}^2}\right) \vphantom{\abs{\hh - \frac{\HH}{n-1}g^\top}^2}\right\rbrace
\end{split}
\end{align}
where $\hh$ and $\HH$ are respectively the second fundamental form and the mean curvature of the level sets of $\varphi$ with respect to the unit normal $\nabla \varphi / \abs{\nabla \varphi}$, $\nabla^\top$ is defined in \cref{def:Geometry_and_orthogonal_decomposition} and $\Ric_g$ denotes the Ricci tensor of the background metric. In particular, $\div(X) \geq 0 $ for $(n-p)/[(n-1)(p-1)]\leq \beta <+\infty$.
\end{lemma}

\proof The proof follows the same lines of \cite[Lemma 4.1]{Agostiniani2019}, replacing accordingly the vector fields $W=\abs{\nabla \varphi}^{p-2} \nabla \abs{\nabla \varphi}^{\beta(p-1)}$ and $Z=(p-2) \abs{ \nabla \varphi}^{p-2} \nabla^\perp \abs{ \nabla \varphi}^{\beta(p-1)}$. The Ricci curvature term appear computing the divergence of $W$ thanks to the Bochner identity for $p$-harmonic functions, as the reader can see following \cite[Proposition 4.3]{Fogagnolo2019}. \endproof

Suppose that $\abs{\nabla \varphi}\neq0$ everywhere. We can apply the Divergence Theorem in the domain $\set{s< \varphi< S}$ to obtain
\begin{equation}\label{eq:difference_derivative_smooth_case}
\int\limits_{\set{\varphi=S}} \kern-.2cm\ip{X, \frac{\nabla \varphi}{\abs{\nabla \varphi}}} \dd \sigma -\kern-.2cm\int\limits_{\set{\varphi=s}} \kern-.2cm\ip{X, \frac{\nabla \varphi}{\abs{\nabla \varphi}}} \dd \sigma =\kern-.2cm \int\limits_{\set{s< \varphi < S}}\kern-.2cm\div X \dd \mu\geq 0.
\end{equation}
Using \eqref{eq:derivativeidentity} we deduce that
\begin{equation}\label{eq:IncreasingDerivativePhi}
\ee^{-\frac{(n-p)}{(n-2)(p-1)}s}(\Phi^\beta_p)'(s)\leq\ee^{-\frac{(n-p)}{(n-2)(p-1)}S}(\Phi^\beta_p)'(S).
\end{equation}
This almost concludes the prove of the Monotonicity Theorem for $\Phi^\beta_p$ with $(n-p)/[(n-1)(p-1)]< \beta <+\infty$ assuming the absence of critical points. Indeed, by integrating it, monotonicity will follow as in \cite[Theorem 3.4]{Fogagnolo2019}. This case lies in the same trail blazed in \cite{Agostiniani2016} since if $\abs{\nabla \varphi} \neq 0$ the $p$-Laplace operator is elliptic non degenerate, and thus the techniques used for harmonic functions fit perfectly.

 If we want to pursue the previous path, even when the critical set of $\phi$ is not empty, we are first committed to provide a version of \eqref{eq:difference_derivative_smooth_case} that holds even in presence of critical values. The main issue is that $\div(X)$ does not belong to $L^1_{\loc}$ a priori. Following the same lines of \cite[Proposition 4.6]{Gigli2021}, testing $s \mapsto\ee^{-\frac{(n-p)}{(n-2)(p-1)}s}(\Phi^\beta_p)'(s)$ against nonnegative functions $\eta \in \CC^\infty_c(0,+\infty)$ and using the coarea formula \cref{prop:coarea} for $f=\ip{X,\nabla \varphi}$ one gets 
\begin{equation}
    (p-1)\int\limits_0^{+\infty} \eta'(s) \ee^{-\frac{(n-p)}{(n-2)(p-1)}s}(\Phi^\beta_p)'(s) \dd s =\kern-.2cm \int\limits_{M\smallsetminus\Crit(\varphi)}\kern-.5cm \ip{X, \nabla [\eta(\varphi)]} \dd \mu.
\end{equation} 
We now would like to integrate by parts and use the nonnegativitiy of $\div(X)$ outside the critical set of $\varphi$. In doing this, we are hampered by the fact that $\div( \chi_{M \smallsetminus \Crit \varphi} X)$ is actually a measure that is possibly not absolutely continuous. Hence we can aim to prove that $s \mapsto\ee^{-\frac{(n-p)}{(n-2)(p-1)}s}(\Phi^\beta_p)'(s)$ belongs to $\mathrm{BV}_{\loc}(0,+\infty)$, but not the absolute continuity. Differently from the nonsmooth case, we can here employ the higher regularity of $\varphi$ outside its critical set to refine the result.

\begin{proposition}\label{prop:IIderivativeidentiy}
Let $(M,g)$ be a $p$-nonparabolic Riemannian manifold with $\Ric\geq 0$. Let $\Omega \subseteq M$ be a open bounded subset with smooth boundary. For every $(n-p)/[(n-1)(p-1)]< \beta <+\infty $, the function $s \mapsto\ee^{-\frac{(n-p)}{(n-2)(p-1)}s}(\Phi^\beta_p)'(s)$ defined in \eqref{eq:derivativeidentity} belongs to $W^{1,1}_{\loc}(0,+\infty)$ and its derivative is given by
\begin{equation}\label{eq:IIderivativeidentiy}
\left(\ee^{-\frac{(n-p)}{(n-2)(p-1)}s}(\Phi^\beta_p)'(s)\right)' =\frac{1}{p-1} \int\limits_{\set{\varphi=s}} \frac{\div X}{\abs{\nabla \varphi}} \dd \sigma,
\end{equation}
for almost every $s \in [0,+\infty)$, where $X$ is the vector field defined in \eqref{eq:Xdef}.
\end{proposition}

\proof \cref{prop:IIderivativeidentiy} follows if we prove that $\div(X)(1 - \chi_{\Crit(\varphi)}) $ belongs to $L^1_{\loc}(M \smallsetminus \overline{\Omega})$ and
\begin{equation}\label{eq:nonnegative_function_duality}
    (p-1)\int\limits_0^{+\infty} \eta'(s) \ee^{-\frac{(n-p)}{(n-2)(p-1)}s}(\Phi^\beta_p)'(s) \dd s  = -\kern-.3cm \int\limits_{M\smallsetminus\Crit\varphi}\kern-.3cm \eta(\varphi)\div X \dd \mu.
\end{equation}
holds for every $\eta\in\CC_c^\infty(0,+\infty)$. Indeed, by coarea formula \eqref{prop:coarea} with $f= \div(X)(1 - \chi_{\Crit(\varphi)})  $ we would get
\begin{equation}
\int\limits_{M\smallsetminus\Crit\varphi}\kern-.3cm \eta(\varphi)\div X \dd \mu= \int\limits_0^{+\infty} \eta(s)\int\limits_{M\smallsetminus\Crit\varphi}\kern-.3cm \div X \dd \sigma \dd t 
\end{equation}
which implies both that $\ee^{-(n-p)s/(n-2)(p-1)}(\Phi^\beta_p)' \in W^{1,1}_{\loc}(0,+\infty)$ and \eqref{eq:IIderivativeidentiy}.

\step[Proof for nonnegative $\eta$] Let $\eta \in \CC_c^\infty(0,+\infty)$ be nonnegative. For any given $\varepsilon>0$ consider the smooth nonnegative cut-off function $\chi_\varepsilon:[0,+\infty)\to \R$ defined as
\begin{equation}\label{eq:cut-off_definition}
\begin{cases}[c@{\ \ }l]
\chi_\varepsilon(t)=0 & \text{in } t< \frac{1}{2}\varepsilon\text{,}\\[.2cm]
0<\chi'_\varepsilon(t)\leq 2\varepsilon^{-1} & \text{in } \frac{1}{2} \varepsilon\leq t \leq \frac{3}{2}\varepsilon,\\[.2cm]
\chi_\varepsilon(t)=1 &\text{in } t>\frac{3}{2}\varepsilon.
\end{cases}
\end{equation}
Consider accordingly the vector field $X_\varepsilon = \chi_\varepsilon( \abs{\nabla \varphi}^{\beta(p-1)}) X$ where $X$ is the vector field given in \eqref{eq:Xdef}. Let $\eta\in \CC^\infty_c(0,+\infty)$ be nonnegative. We notice that $\abs{ \ip{X_\varepsilon, \nabla \varphi}} \leq \abs{ \ip{X, \nabla \varphi}}$ which is in $L^2_{\loc}(M \smallsetminus \overline{\Omega})$ by \eqref{gradphi-bound} and \cref{rmk:vanishing_of_gradient}. Hence \eqref{eq:derivativeidentity}, coarea formula with $f= \eta'(\varphi) \ip{X,\nabla \varphi/ \abs{\nabla \varphi}}$ and Dominated Convergence Theorem imply
\begin{equation}
\displaystyle\int\limits_0^{+\infty}\eta'(s) \ee^{-\frac{(n-p)}{(n-2)(p-1)}s}(\Phi^\beta_p)'(s) \dd s =\lim_{\varepsilon\to 0^+}\frac{1}{p-1}\displaystyle\displaystyle\int\limits_M\eta'(\varphi) \ip{X_{\varepsilon}, {\nabla \varphi}} \dd \mu.
\end{equation}
Employing the coarea formula in \eqref{eq:enached_coarea} with $f = \ip{X_\varepsilon, \nabla \varphi}$ and integration by parts, we obtain that
\begin{align}
&\int\limits_0^{+\infty} \eta'(s)\kern-.2cm\int\limits_{\set{\varphi=s}}\kern-.2cm\ip{ X_{\varepsilon}, \frac{\nabla \varphi}{\abs{\nabla \varphi}}} \dd \sigma\dd s = \int\limits_{M} \eta'(\varphi) \ip{X_ \varepsilon, \nabla \varphi} \dd \mu= -\int\limits_{M}\div( X_\varepsilon)\eta(\varphi) \dd \mu\\
&\kern.5cm=-\kern-.3cm\int\limits_{M \smallsetminus N_{\varepsilon/2}}\kern-.3cm\eta(\varphi) \chi_\varepsilon(\abs{\nabla \varphi}^{\beta(p-1)}) \div X \dd \mu -\kern-.5cm \int\limits_{N_{3 \varepsilon/2} \smallsetminus N_{\varepsilon/2}}\kern-.5cm\eta(\varphi) \chi'_\varepsilon(\abs{\nabla \varphi}^{\beta(p-1)})\ip{X, \nabla \abs{\nabla \varphi}^{\beta(p-1)}} \dd \mu
\end{align}
where $N_\delta = \set{ \abs{\nabla \varphi}^{\beta (p-1)} < \delta}$ for every $\delta >0$. By Monotone Convergence Theorem, the first integral gives
\begin{equation}
\lim_{\varepsilon\to 0^+}\kern-.2cm\int\limits_{M \smallsetminus N_{\varepsilon/2}}\kern-.3cm\eta(\varphi) \chi_\varepsilon(\abs{\nabla \varphi}^{\beta(p-1)}) \div X \dd \mu=\kern-.5cm \int\limits_{M\smallsetminus \Crit(\varphi)}\kern-.5cm\eta(\varphi) \div X \dd \mu\geq 0
\end{equation}

To conclude the step, it thus remains to prove that the second integral vanishes as $\varepsilon\to 0^+$. Firstly, we observe that $\abs{\nabla \varphi}^{\beta (p-1)}\geq \varepsilon/2$ on $N_{3 \varepsilon/2} \smallsetminus N_{\varepsilon/2}$, then $\varphi$ is smooth. By $\ip{X, \nabla \abs{\nabla \varphi}^{\beta(p-1)}}\geq 0$ and coarea formula in \cref{prop:coarea} we obtain that
\begin{align}
&\abs{\kern.1cm\int\limits^{}_{N_{3 \varepsilon/2} \smallsetminus N_{\varepsilon/2}}\kern-.5cm \eta(\varphi)\chi'_\varepsilon(\abs{\nabla \varphi}^{\beta(p-1)})\ip{X, \nabla \abs{\nabla\varphi}^{\beta(p-1)}} \dd \mu\kern.1cm}\\&\kern6cm \leq  \frac{2}{\varepsilon}\norm{\eta}_{L^\infty} \int\limits_{\varepsilon/2}^{3\varepsilon/2}\int\limits_{\partial N_s}\frac{\displaystyle\ip{X, \nabla \abs{\nabla \varphi}^{\beta(p-1)}}}{\abs{\nabla \abs{\nabla \varphi}^{\beta(p-1)}}} \dd \sigma \dd s.
\end{align}
Let $R>0$ and $\mathcal{H}$ be defined as
\begin{equation}
\mathcal{H}(r)= \int\limits_{\partial N_r} \frac{\displaystyle\ip{X, \nabla \abs{\nabla \varphi}^{\beta(p-1)}}}{\abs{\nabla \abs{\nabla \varphi}^{\beta(p-1)}}} \dd \sigma
\end{equation}
for every $r \in (0,R)$ regular value of $\abs{\nabla \varphi}$, hence for almost every $r\in(0,R)$ thanks to Sard's Theorem. By the Mean Value Theorem, showing that $\mathcal{H}(r)$ vanishes as $r \to 0^+$ is enough to conclude the proof. 

Let $0<t< r<R$ be two regular values for $\abs{\nabla \varphi}$, applying the Divergence Theorem to the smooth vector field $X$ on $N_r \smallsetminus N_t$ we get
\begin{align}
\mathcal{H}(r)-\mathcal{H}(t)&= \int\limits_{\partial N_r} \frac{\displaystyle\ip{X, \nabla \abs{\nabla \varphi}^{\beta(p-1)}}}{\abs{\nabla \abs{\nabla \varphi}^{\beta(p-1)}}} \dd \sigma - \int\limits_{\partial N_t} \frac{\displaystyle\ip{X, \nabla \abs{\nabla \varphi}^{\beta(p-1)}}}{\abs{\nabla \abs{\nabla \varphi}^{\beta(p-1)}}} \dd \sigma \\&= \int\limits_{N_r \smallsetminus N_t} \div(X) \dd \mu
=\int\limits_t^r\int\limits_{\partial N_s}\frac{\div(X)}{\abs{\nabla\abs{\nabla \varphi}^{\beta(p-1)}}} \dd \sigma \dd s.\label{eq:FasIntegral}
\end{align}
where the last identity is an application of the coarea formula. Since the integrand in the rightmost side is nonnegative and $\mathcal{H}$ is almost everywhere finite, $\mathcal{H}$ is locally absolutely continuous. 

Since $\Ric\geq 0$ and $\abs{\nabla \varphi}^2\abs{\hh - \frac{\HH}{n-1}g^\top}^2 \geq 0$, by \eqref{eq:Qdef} we have that
\begin{align}
\div X & \geq \ee^{-\frac{(n-p)}{(p-1)(n-2)}S} \beta(p-1)\abs{\nabla \varphi}^{\beta(p-1)+p-4  }\left((p-1)\left[\beta+\frac{p-2}{p-1}\right]\right.\\&\kern2cm \abs{\nabla^\top\abs{\nabla \varphi}}^2 
\left.+(p-1)^2 \left[\beta - \frac{(n-p)}{(p-1)(n-1)}\right]\abs{\nabla^\perp\abs{\nabla \varphi}}^2\right)\\
 & \geq  \beta^2 (p-1)^2\Cn \abs{ \nabla \varphi}^{\beta(p-1)+ p -4}\left(\abs{\nabla^\perp\abs{\nabla \varphi}}^2+\abs{\nabla^\top\abs{\nabla \varphi}}^2 \right)\\
 &\geq \Cn  \abs{\nabla \varphi}^{-\beta(p-1)+p-2}\abs{ \nabla \abs{ \nabla \varphi}^{\beta(p-1)}}^2 
\end{align}
where
\begin{align}
\Cn=\frac{1}{\beta} \ee^{-\frac{(n-p)}{(p-1)(n-2)}S} \min\left\lbrace\left[ \beta + \frac{p-2}{p-1}\right], (p-1)\left[\beta- \frac{(n-p)}{(p-1)(n-1)}\right] \right\rbrace>0
\end{align}
 Taking derivatives in \eqref{eq:FasIntegral} it holds that
\begin{equation}
\mathcal{H}'(r)=\int\limits_{\partial N_s}\frac{\div(X)}{\abs{\nabla\abs{\nabla \varphi}^{\beta(p-1)}}} \dd \sigma \geq  \Cn\int\limits_{\partial N_r}  \abs{\nabla \varphi}^{-\beta(p-1)+p-2}\abs{\nabla \abs{\nabla \varphi}^{\beta(p-1)}} \dd \sigma = \Cn\, \frac{\mathcal{H}(r)}{r}
\end{equation}
for almost any $r>0$. Integrating for $R>r$, it holds that
\begin{equation}
\frac{\mathcal{H}(r)}{r^{\Cn}} \leq \frac{\mathcal{\mathcal{H}}(R)}{R^{\Cn}}.
\end{equation}
We thus deduce that $\mathcal{H}(r)\to 0$ as $r \to 0^+$. 

\step[Conclusions] In the previous step we proved \eqref{eq:nonnegative_function_duality} for every nonnegative function $\eta\in \CC^\infty_c(0,+\infty)$. Let be $K\subset M \smallsetminus \overline{\Omega}$. Then, there exists a $\eta_K \in \CC^\infty_c(0,+\infty)$, $\eta_K\geq 0$, such that $\eta_K(\varphi)\geq 1 $ on $K$. In particular, since $\div(X)\geq 0$ outside $\Crit(\varphi)$ we have
\begin{align}
    \int\limits_{K} \div(X) (1 - \chi_{\Crit(\varphi)}) \dd \mu &\leq \int\limits_{M \smallsetminus \Crit(\varphi)}\kern-.5cm \eta_K(\varphi) \div(X) \dd \mu\\
    &= -(p-1)\int\limits_0^{+\infty} \eta_K'(s) \ee^{-\frac{(n-p)}{(n-2)(p-1)}s}(\Phi^\beta_p)'(s) \dd s
\end{align}
which is finite thanks to \cref{prop:derivativeidentity}. This ensures that $\div(X)(1 - \chi_{\Crit(\varphi)})$ belongs to $L^1_{\loc}(M \smallsetminus \overline{\Omega})$. Approximating the positive and the negative part of a general $\eta\in \CC^\infty_c(0,+\infty)$, that are nonnegative Lipschitz with compact support, we can conclude. \endproof

\proof[Proof of \cref{prop:Monotonicity_Fbeta_sequence}] We use an argument due to Colding and Minicozzi in \cite{Colding2014}. By \cref{prop:IIderivativeidentiy,prop:derivativeidentity}, $\Phi_p^\beta$ is $W^{2,1}_{\mathrm{loc}}(0, + \infty)$. By \eqref{eq:IIderivativeidentiy} $s \mapsto\ee^{-\frac{(n-p)}{(n-2)(p-1)}s}(\Phi^\beta_p)'(s)$ is nondecreasing, then for every $0\leq s<S<+\infty$ we have 
\begin{equation}
\ee^{-\frac{(n-p)}{(n-2)(p-1)}(S-s)}(\Phi^\beta_p)'(s)\leq (\Phi^\beta_p)'(S).
\end{equation}
Integrating the above inequality, we get
\begin{equation}\label{eq:derivativeinequality}
\frac{(n-1)(p-1)}{(n-p)}\left(\ee^{\frac{(n-p)}{(n-2)(p-1)}(S-s)}-1\right)(\Phi_p^\beta)'(s)\leq \Phi_p^\beta(S)-\Phi_p^\beta(s)
\end{equation}
for every $0\leq s<S<+\infty$. Suppose, by contradiction, that $(\Phi_p^\beta)'(s)>0$ for some $s\in [0,+\infty)$. Passing to the limit as $S\to +\infty $ in \eqref{eq:derivativeinequality} we would get that $\Phi_p^\beta(S)\to+\infty$ against the boundedness property ensured by \cref{bound-Phi}. Hence, $(\Phi_p^\beta)'(s)\leq 0$ and in particular $s \mapsto \Phi_p^\beta(s)$ is nonincreasing. Notice that $\Phi^\beta_p$ is a bounded, nonincreasing $\CC^1(0,+\infty)$ function, then $(\Phi^\beta_p)'(s)\to 0$ as $s \to +\infty$. Coupling this with the coarea formula in \cref{prop:coarea} for $f= \div(X)(1- \chi_{\Crit(\varphi)})/ \abs{\nabla \varphi} $ one get that
\begin{equation}\label{eq:last_integral_identiy}
\begin{split}
    \ee^{-\frac{(n-p)}{(n-2)(p-1)}s}(\Phi^\beta_p)'(s) &= \lim_{S\to +\infty} \ee^{-\frac{(n-p)}{(n-2)(p-1)}s}(\Phi^\beta_p)'(s) -\ee^{-\frac{(n-p)}{(n-2)(p-1)}S}(\Phi^\beta_p)'(S) \\&=\lim_{S \to +\infty} 
    \kern.2cm-\kern-1cm\int\limits_{\set{s\leq\varphi \leq S}\smallsetminus \Crit(\varphi)}\kern-.8cm \div X \dd \mu=\kern.3cm-\kern-.7cm\int\limits_{\set{\varphi \geq s}\smallsetminus \Crit(\varphi)}\kern-.6cm \div X \dd \mu,
\end{split}
\end{equation}
which also ensures that $\div X \in L^1(M\smallsetminus (\overline{\Omega} \cup \Crit (\varphi)))$.

One can now obtain \eqref{eq:derivative_of_Fbeta} rewriting \eqref{eq:last_integral_identiy} in terms of $u$. The proof proceeds through direct computations. The main ones are contained in \cite[Section 3.3]{Fogagnolo2019}, the only difference is the Ricci term that can be computed as
\begin{equation}
\Ric(\nabla \varphi, \nabla \varphi)= \left[\frac{(p-1)(n-2)}{(n-p)}\right]^2 u^{-2 \frac{n+p-2}{n-p}} \Ric(\D u, \D u).
\end{equation}
Consequently, \eqref{eq:2nd_derivative_of_Fbeta} follows by \eqref{eq:derivative_of_Fbeta} and coarea formula.

For the rigidity statement, suppose that $(F_p^\beta)'(t_0)=0$ for some $t_0 \in [1,+\infty)$ regular for $u$. Then by \eqref{eq:derivative_of_Fbeta}
\begin{align}
\abs{ \hh- \frac{\HH}{n-1}g^{\top}}_g=0, &&\abs{ \D^\top \abs{ \D u}_g}_g=0
\end{align}
hold on $\set{u \leq 1/t_0}\smallsetminus \Crit(u)$. By \cref{prop:Kato-type}, $(\set{u\leq 1/t_0},g)$ splits to a warped product near the level set $\set{u = 1/t_0}$. In particular, the mean curvature $\HH$ depends only on $u$. By \eqref{eq:derivative_of_F_standard} also $\abs{ \D u } $ depends only on $u$ and
\begin{equation}
    \frac{ \partial }{\partial u} \abs{ \D u}_g= \frac{ \HH}{p-1}= \frac{n-1}{n-p} \frac{\abs{\D u}_g}{u}.
\end{equation}
Integrating it we get that for some $A(t_0)>0$ the identity 
\begin{equation}
    \abs{ \D u}_g = u^{\frac{n-1}{p-1}}A(t_0)
\end{equation}
holds, which gives that $\abs{\D u}_g$ never vanishes on $\set{u \leq 1/t_0}$ by the continuity of gradient. Recalling the relation between $u$, $\eta$ and $t$ in \eqref{eq:Kato-type_rigidity}, we obtain that $\eta(t)= B(t_0)t_0t + (1-B(t_0))$ for some $B(t_0)>0$. If we define the new coordinate as $\tau = t +[1-B(t_0)]/[B(t_0)t_0]$ and $\tau_0 =1/[t_0(B(t_0)-1)]$, we have that $\set{\tau\geq \tau_0}= \set{ u \leq 1/t_0}$, $\eta(t)=\tau/\tau_0$ and $\dd \tau = - \dd t$. To sum up, we have proven that $(\set{u\leq 1/t_0},g)$ is isometric to 
\begin{equation}
    \left([\tau_0, +\infty) \times \set{u = 1/t_0}, \dd \tau \otimes \dd \tau + \left(\frac{\tau}{\tau_0}\right)^2 g_{\set{u=1/t_0}} \right)
\end{equation}
leaving us only to characterise $\tau_0$.
Observe that, by the conical splitting, the measure of the level sets of $\tau$ satisfy
\begin{equation}
\abs{\set{\tau=R}}=\left(\frac{R}{\tau_0}\right)^{n-1} \abs{\set{u=1/t_0}}.
\end{equation}
One can easily prove that on a cone
\begin{equation}
1=\lim_{R\to +\infty}\frac{\abs{\set{\tau \leq R}}}{\abs{B(o,R)}}=\lim_{R\to +\infty}\frac{\abs{\set{\tau= R}}}{\abs{\partial B(o,R)}}
\end{equation}
that can be used to compute the claimed value of $\tau_0$
\begin{equation}
\AVR(g) = \lim_{R \to+ \infty} \frac{\abs{\set{\tau=R}}}{R^{n-1} \abs{\S^{n-1}}}= \frac{\abs{\set{u=1/t_0}}}{\tau_0^{n-1} \abs{\S^{n-1}}}.
\end{equation} 
\endproof

We conclude this section sketching the proof of the Monotonicity-Rigidity Theorem for $\Phi^\infty_p$, which does not require much more effort than in $\R^n$ \cite{Fogagnolo2019}.

\proof[Proof of \cref{prop:Monotonicity_Finfty_sequence}] Lemma 5.1 in \cite{Fogagnolo2019} holds also in this setting. The only difference in proving that $\abs{\nabla \varphi}^p$ is a subsolution of the nondegenerate uniformly elliptic operator 
\begin{equation}
\mathscr{L}(f)= \Delta f +(p-2) \nabla \nabla f \left( \frac{ \nabla \varphi}{\abs{\nabla \varphi}}, \frac{ \nabla \varphi}{\abs{\nabla  \varphi}}\right) - \frac{n-p}{n-2} \ip{\nabla f, \nabla \varphi},
\end{equation}
acting on smooth $f$ in a neighbourhood of points such that $\abs{\nabla \varphi}>0$, is that the curvature term that appear when Bochner identity for $p$-harmonic functions is applied can be controlled by $\Ric \geq 0$. We claim that
\begin{equation}\label{eq:claim_Phiinfty}
\abs{\nabla \varphi}(x)\leq\sup_{ \set{ \varphi = s}}\abs{\nabla \varphi}
\end{equation}
for every $s\in[0,+\infty)$ and $x \in \set{\varphi \geq s }$, which is the main ingredient in the proof of \cite[Theorem 3.5]{Fogagnolo2019}. Firstly suppose that $\Phi^\infty_p(s)>0$ and let be $0<\delta <\Phi^\infty_p(s)$. By \eqref{bound-Phi} $\abs{\nabla \varphi}\leq \Cn$ uniformly in $M\smallsetminus \Omega$. For some $S>s$ let 
\begin{equation}
w= \abs{\nabla \varphi}^p - \sup_{\set{\varphi =s}} \abs{\nabla \varphi}^p - \Cn^p \ee^{\frac{n-p}{(n-2)(p-1)}(\varphi -S)}
\end{equation}
be defined on $\set{s \leq \varphi \leq S}\smallsetminus N_\delta$ where $N_\delta = \set{\abs{\nabla \varphi}<\delta}$. Since $w \leq 0$ on the boundary of $\set{s \leq \varphi \leq S}\smallsetminus N_\delta$ and $\mathscr{L}(w)\geq 0$ in its interior, by the Maximum Principle we have that
\begin{equation}\label{eq:auxiliary_maximum_principle}
\abs{\nabla \varphi}^p\leq \sup_{\set{\varphi =s}} \abs{\nabla \varphi}^p + \Cn^p \ee^{\frac{n-p}{(n-2)(p-1)}(\varphi -S)}
\end{equation} 
on $\set{s \leq \varphi \leq S} \smallsetminus N_\delta$. Moreover, since $\abs{\nabla \varphi}< \delta $ on $N_\delta$, \eqref{eq:auxiliary_maximum_principle} is thus satisfied in the whole $\set{s \leq \varphi \leq S}$. Passing to the limit as $S \to +\infty$, \eqref{eq:claim_Phiinfty} is proven for $s\in[0,+\infty)$ such that $\Phi^\infty_p(s)>0$.

We now prove \cref{cor:regualirity_of_p-harmonic}, namely that $\Phi^\beta_p(s)>0$ for every $s \in [0,+\infty)$ that in particular yields \eqref{eq:claim_Phiinfty} proving the monotonicity. Suppose by contradiction that $\Phi^\infty_p(s)=0$ for some $s \in [0,+\infty)$. By \cref{prop:coarea} there exists a sequence of $(s_j)_{j\in \N}$, $s_j \to s$ as $j\to +\infty$ and $\Phi^\infty_p(s_j)>0$. If, up to a subsequence, we can assume that $\Phi^\infty_p(s_j)\to 0$, then we can conclude. Indeed, $\Phi^\infty_p(s_j)\geq \abs{\nabla \varphi}(x)$ for every $x \in \set{\varphi \geq s}$ and $\Phi^\infty_p(s_j)\to 0$ as $j\to +\infty$, hence $\abs{\nabla \varphi}=0$ on $\set{\varphi \geq s}$, contradicting the unboundedness of $\varphi$. Suppose now that every subsequence of $\Phi^\infty_p(s_j)$ does not vanish, then there would be a $\delta>0$ and $J\in \N$ such that $\Phi_p^\infty(s_j)>\delta$ for every $j\geq J$. Since level sets of $\phi$ are compact, $\Phi_p^\beta(s_j)$ is actually achieved at some point $x_{s_j}\in \set{ \varphi =s_j}$. Moreover, $(x_{s_j})_{j \in \N}$ is bounded, since it is contained in $\set{\varphi \leq s}$. Hence, we can assume that there exists $x \in \set{\varphi \leq s}$ such that $x_{s_j}\to x$ as $j \to +\infty$. Since $\varphi$ is $\CC^1$, we obtain that $\varphi(x)=s$ and $ \abs{ \nabla \varphi}(x)\geq\delta$, contradicting the fact that $\Phi^\infty_p(s)=0$. 

Using a similar argument we can infer that $s\mapsto \Phi_p^\beta(s)$ is left continuous. Indeed, by contradiction there would be a $\delta>0$ such that $\Phi_p^\infty(s)\geq \Phi_p^\infty(s_0)+ \delta$ for any $s<s_0$. Let $x_s\in \set{\varphi=s}$ such that $\Phi^\infty_p(s)= \abs{\nabla \varphi}(x_s)$. By the compactness of $\set{ \varphi \leq s_0}$, there exists a sequence $(s_j)_{j \in \N}$ and a point $x \in \set{\varphi \leq s_0}$ such that $s_j <s_0$, $s_j \to s_0$ and $x_{s_j} \to x$. Since $\varphi \in \CC^1$, $\varphi(x)= s_0$ and $\abs{\nabla \varphi}(x)\geq\Phi_p^\infty(s_0)+\delta$, contradicting the definition of $\Phi_p^\infty$. To prove the right continuity it is the enough to prove that $s \mapsto \Phi_p^\infty(s)$ is lower semicontinuous. Since $\Phi_p^\infty>0$, the maximum of $\abs{ \nabla \varphi}$ on $\set{\varphi=s} $ is achieved at a regular point $x$. Let $(s_j)_{j \in \N}$ be a sequence such that $s_j \to s$ as $j \to +\infty$. Seeing as $\abs{\nabla \varphi}$ is continuous, there exists a sequence of points $(x_{s_j})_{j \in \N}$ such that $x_{s_j}\in \set{\varphi=s_j}$ and $x_{s_j} \to x$ as $j \to +\infty$. Since $\abs{ \nabla \varphi}(x_{s_j}) \leq \Phi^\infty_p (s_j)$ for every $j \in \N$, we complete the proof.

We turn to prove the second part of \cref{prop:Monotonicity_Finfty_sequence}. Since $x_t$ is a point of maximum for the function $\abs{\D u }_g/u^{(n-1)/(n-p)}$ on $\set{u \leq 1/t}$, its derivative with respect to the normal unit vector $\nu_t= - \D u / \abs{ \D u}_g$ is nonpositive. Hence \eqref{eq:mean_curvature_inequality_background} follows by direct computations. To conclude, both rigidity statements follow in the same way as in \cite[Theorem 3.5]{Fogagnolo2019}, since $\abs{\D u }_g^p/u^{p(n-1)/(n-p)}$ is also a subsolution of $\mathscr{L}f=0$, thanks to \eqref{gradphi-relation}. \endproof

\section{Geometric consequences of the Monotonicity Theorems}
\label{sec:consequences}
In this section we prove the geometric implications of the Monotonicity-Rigidity theorems, that are the Minkowski Inequalities, a rigidity result under pinching condition and a sphere theorem. The proof of these theorems follow, along with the Monotonicity already mentioned, by a contradiction argument that involves the Iso-$p$-capacitary Inequality, that we are going to state and prove immediately, since we believe of independent interest.

\subsection{Iso-\texorpdfstring{$p$}{p}-capacitary Inequality}
We provide the sharp Iso-$p$-capacitary Inequality in complete noncompact Riemannian manifolds with nonnegative Ricci curvature and Euclidean volume growth. As for the standard Iso-$p$-capacitary Inequality in Euclidean setting, the proof fully relies on the Isoperimetric Inequality combined with a P\'olya-Szeg\"o principle. In particular, the sharpness of the inequality that follows is a direct consequence of the sharp isoperimetric constant in this setting, that has been found first in dimension $3$ in \cite{Agostiniani2018} and later extended to all dimensions in \cite{Brendle2020}. See also \cite{Fogagnolo2020a, Balogh2021, Johne2021} for related results. The proof below is classical, and it is inspired by \cite{Jauregui2012}, where it is illustrated for the $2$-capacity in $\R^n$.

\begin{theorem}[Iso-$p$-capacitary Inequality]\label{thm:isopcapacitary}
Let $(M,g)$ be a complete, noncompact Riemmannian manifold with nonnegative Ricci curvature and Euclidean Volume Growth. Let be $\Omega\subseteq M$ open bounded subset with smooth boundary. Then
\begin{equation}\label{eq:isopcap}
   \frac{\capa_p(\mathbb{B}^n)^n}{\abs{\mathbb{B}^n}^{n-p}}\AVR(g)^{p}\leq\frac{\capa_p(\Omega)^n}{\abs{\Omega}^{n-p}}
\end{equation}
Moreover, if the equality holds then $(M,g)$ is isometric to the Euclidean Space and $\Omega$ is a ball.
\end{theorem}

\proof By \eqref{cap-u} and the coarea formula in \cref{prop:coarea} we have that
\begin{equation}\label{eq:coarea_formula_capacity}
    \capa_p(\Omega)= \int\limits_{M \smallsetminus \overline{\Omega}} \abs{\D u}^p \dd \mu=\int\limits_0^1 \int\limits_{\set{u=\tau}}\abs{\D u}^{p-1} \dd \sigma \dd \tau.
\end{equation}
The H\"{o}lder inequality with exponents $a=p$ and $b=p/(p-1)$ gives
\begin{equation}\label{eq:level_area_bound_CS}
    \abs{ \set{u=\tau}}^p \leq \left(\kern.1cm\int\limits_{ \set{u=\tau}} \abs{\D u}^{p-1} \dd \sigma \right)\left(\kern.1cm \int\limits_{\set{u=\tau}}\frac{1}{\abs{ \D u }}\dd \sigma\right)^{p-1}
\end{equation}
for almost every $\tau \in (0,1]$. Let $V':(0,1]\to \R$ be defined as
\begin{equation}\label{eq:derivative_volume_superlevel}
    V'(\tau)= -\int\limits_{\set{u=\tau}} \frac{1}{\abs{\D u}} \dd \sigma
\end{equation}
Moreover, let $V:(0,1]\to \R$ be the primitive of $V'(\tau)$ chosen as
\begin{equation}\label{eq:VolumeIsopDefinition}
    V(\tau)= \abs{\Omega}-\int\limits_\tau^1 V'(s) \dd s= \abs{ \Omega_\tau \smallsetminus \Crit(u)},
\end{equation}
where the second identity is obtained coupling \eqref{eq:derivative_volume_superlevel} with the coarea formula \eqref{eq:enached_coarea} applied with $f= (1-\chi_{\Crit(u)}) \abs{ \D u}^{-1}$  (see \cref{rmk:coarea_for_fraction}). 

By the Isoperimetric Inequality in \cite[Corollary 1.3]{Brendle2020} we have that
\begin{equation}\label{eq:isopforisopcap}
        \abs{ \set{u=\tau}} \geq \abs{\partial \Omega_\tau} \geq \abs{ \Omega_\tau}^{\frac{n-1}{n}} \AVR(g)^{\frac{1}{n}} n \abs{\mathbb{B}^n}^{\frac{1}{n}}\geq V(\tau)^{\frac{n-1}{n}}\AVR(g)^{\frac{1}{n}} n \abs{\mathbb{B}^n}^{\frac{1}{n}}.
\end{equation}
Let $R(\tau)$ be the radius of the ball in $\R^n$ which has volume $V(\tau)$, then $V(\tau)= \abs{ \mathbb{B}^n} R(\tau)^n$ and $V'(\tau)= \abs{ \Sf^{n-1}}R(\tau)^{n-1} R'(\tau)$. Coupling \eqref{eq:isopforisopcap} with \eqref{eq:coarea_formula_capacity}, \eqref{eq:level_area_bound_CS} and \eqref{eq:derivative_volume_superlevel} we obtain
\begin{align}
    \capa_p(\Omega)&\geq \int\limits_0^{1}\frac{\abs{\set{u=\tau}}^p}{[-V'(\tau)]^{p-1}} \dd \tau\geq n^p\left( \abs{\mathbb{B}^n}\AVR(g)\right)^{\frac{p}{n}}  \int\limits_0^1 \frac{V(\tau)^{\frac{p(n-1)}{n}}}{[-V'(\tau)]^{p-1}} \dd \tau\\
    &= \abs{\Sf^{n-1}}\AVR(g)^{\frac{p}{n}}\int\limits_0^1 \frac{R(\tau)^{n-1}}{[-R'(\tau)]^{p-1}} \dd \tau.
\end{align}
Let now  $v: \set{\abs{x}\geq R(1)}\subset \R^n\to (0,1]$ be the function which is $\tau$ on $\set{\abs{x}=R(\tau)}$. By \eqref{eq:isopforisopcap} and \eqref{gradient-bounduf} there exists a positive constant $\Cn= \Cn(p,n)$ such that
\begin{equation}
    -V'(\tau)= \int\limits_{\set{u=\tau}} \frac{1}{\abs{\D u }} \dd \sigma \geq \Cn \abs{\Omega}^{\frac{n-1}{n}}\tau^{\frac{n-p}{p-1}}.
\end{equation}
Seeing as $\abs{\D v}=-1/R'(\tau)=-\abs{\S^{n-1}}R^{n-1}(\tau)/V'(\tau)$ the function $v$ is locally Lipschitz. Since $\abs{\S^{n-1}}R(\tau)^{n-1}=\abs{\set{\abs{x}=R(\tau)}}= \abs{\set{v =\tau}}$ by the coarea formula \eqref{eq:enached_coarea} applied with $f=\abs{ \D v}^{p-1}$(see \cref{rmk:coarea_for_fraction}) we have
\begin{align}
    \abs{\Sf^{n-1}}\AVR(g)^{\frac{p}{n}}&\int\limits_0^1 \frac{R(\tau)^{n-1}}{[-R'(\tau)]^{p-1}} \dd \tau=\AVR(g)^{\frac{p}{n}}\int\limits_0^1 \int\limits_{\set{v=\tau}}\kern-.2cm \abs{\D v}^{p-1} \dd \sigma \dd \tau \\&= \AVR(g)^{\frac{p}{n}}\kern-.4cm\int\limits_{\set{\abs{x}\geq R(1)}}\kern-.4cm \abs{\D v}^p \dd x 
    \geq \AVR(g)^{\frac{p}{n}}\capa_p\left(\set{\abs{x} < R(1)}\right)
\end{align}
the last one is by the definition of the $p$-capacity \eqref{cap} in flat $\R^n$. Using \eqref{eq:scaling-invariant-cap} and the fact that $\abs{\set{\abs{x}\leq R(1)}}=V(1)= \abs{ \Omega}$, we finally obtain
\begin{equation}
    \AVR(g)^{\frac{p}{n}}\capa_p\left(\set{\abs{x} < R(1)}\right)=\AVR(g)^{\frac{p}{n}}\capa_p\left(\mathbb{B}^n\right)R(1)^{n-p}=\AVR(g)^{\frac{p}{n}}\frac{\capa_p(\mathbb{B}^n)}{\abs{\mathbb{B}^n}^{\frac{n-p}{n}}} \abs{ \Omega}^{\frac{n-p}{n}},
\end{equation}
and consequently \eqref{eq:isopcap}.

Clearly, if the equality holds in \eqref{eq:isopcap} then also the equality holds in the use of the Isoperimetric Inequality and \cite[Theorem 1.2]{Brendle2020} forces the rigidity both of the ambient manifold and $\Omega$. \endproof

We conclude this subsection with the following remark, whose importance will be clarified in the very proof of the $L^p$-Minkowski Inequality (\cref{thm:Lp_mink} below), where a sharp lower bound for the $p$-capacity of the super-level sets of the $p$-capacitary potential of $\Omega$ will be needed. 

\begin{remark} We observe that, replacing $\Omega$ and $u$ with $\Omega_t=\set{u>1/t}\cup \Omega$ and $u_t=tu$ respectively and defining $V:(0,1]\to \R$ in \eqref{eq:VolumeIsopDefinition} as
\begin{equation}
        V(\tau)= \abs{\Omega_t\cup\set{u_t=1}}+\int\limits_{\tau}^{1} \int\limits_{\set{u_t=s}} \frac{1}{\abs{\D u_t}} \dd \sigma \dd s= \abs{\Omega_{\tau/t} \smallsetminus (\Crit(u)\cap \set{\tau<u_t<1}) }
\end{equation}
we obtain that
\begin{equation}\label{eq:isopcap_tlevel}
   \frac{\capa_p(\mathbb{B}^n)^n}{\abs{\mathbb{B}^n}^{n-p}}\AVR(g)^{p}\leq\frac{\capa_p(\Omega_t)^n}{\abs{\Omega_t}^{n-p}}
\end{equation}
holds for every $t \in [1,+\infty)$.
\end{remark}

\subsection{Minkowski Inequality}
We are now ready to prove the $L^p$-Minkowski Inequality in our setting. Let $(M,g)$ be a noncompact, complete Riemannian manifold with $\Ric\geq0$ and Euclidean Volume Growth. Consider the function $t\mapsto F_p(t)$ defined in \eqref{eq:monotonefunction} as $F_p^\beta$ with $\beta=1/(p-1)$. By \eqref{eq:scaling-invariant-cap} we can rewrite $F_p$ in a more geometric fashion as
\begin{equation}\label{eq:newgeometricitnerpretation}
    F_p(t)= t^{\frac{n-1}{n-p}} \int\limits_{\set{u=1/t}} \abs{\D u}^p\dd \sigma = \left(\frac{\Cn_p(\Omega_t)}{\Cn_p(\Omega)}\right)^{-\frac{n-p-1}{n-p}} \int\limits_{\set{u_t=1}}\abs{\D u_t}^p\dd \sigma.
\end{equation}
where $u_t= tu$ and $\Omega_t= \set{ u> 1/t}\cup \Omega$.

\begin{theorem}[$L^p$-Minkowski Inequality] \label{thm:Lp_mink}
Let $(M,g)$ be complete Riemannian manifold with $\Ric \geq 0$ and Euclidean Volume Growth. Let $\Omega \subseteq M$ be a open bounded subset with smooth boundary. Then for every $1<p<n$, the following inequality holds
\begin{equation}\label{eq:Lp_minkowski}
\mathrm{C}_p(\Omega)^{\frac{n-p-1}{n-p}}\AVR(g)^{\frac{1}{n-p}}\leq \frac{1}{\abs{\S^{n-1}}} \int\limits_{\partial\Omega} \abs{\frac{\HH}{n-1}}^p \dd \sigma.
\end{equation}
Moreover, the equality holds in \eqref{eq:Lp_minkowski} if and only if $(M\smallsetminus \Omega, g)$ is isometric to 
\begin{align}
    \left( [\rho_0, +\infty) \times \partial \Omega, \dd \rho \otimes \dd \rho + \left( \frac{\rho}{\rho_0}\right)^2 g^{}_{\partial \Omega}\right),&&\text{where }\rho_0 =\left(\frac{\abs{\partial \Omega}}{\AVR(g)\abs{\S^{n-1}}}\right)^{\frac{1}{n-1}}.
\end{align}
\end{theorem}

\proof 
We first show that
\begin{equation}\label{eq:Minkowski-contradicition}
    \mathrm{C}_p(\Omega)^{\frac{n-p-1}{n-p}}\AVR(g)^{\frac{1}{n-p}}\leq \frac{1}{\abs{\S^{n-1}}} \left(\frac{p-1}{n-p}\right)^p \int\limits_{\partial \Omega} \abs{ \D u }^p \dd \sigma
\end{equation}
holds for any open subset $\Omega\subseteq M$ with smooth boundary.

Let then $\otheta < \AVR(g)$ and suppose by contradiction that there exists an open subset $\Omega\subseteq M$ with smooth boundary, such that 
\begin{equation}
\mathrm{C}_p(\Omega)^{\frac{n-p-1}{n-p}}\otheta^{\frac{1}{n-p}}\geq \frac{1}{\abs{\S^{n-1}}} \left(\frac{p-1}{n-p}\right)^p \int\limits_{\partial \Omega} \abs{ \D u }^p \dd \sigma.
\end{equation}
Define $\tau=1/t\in (0,1]$. By \cref{prop:Monotonicity_Fbeta_sequence}, the function $\tau\mapsto F_p(\tau)$ is nondecreasing for $\tau\in(0,1]$. Exploiting this monotonicity as in \eqref{eq:newgeometricitnerpretation} we have
\begin{equation}
    \left(\frac{n-p}{p-1}\right)^p\abs{\Sf^{n-1}}\otheta^{\frac{1}{n-p}}\geq\Cn_p(\Omega)^{-\frac{n-p-1}{n-p}}\int\limits_{\partial \Omega}\abs{\D u}^p \dd \sigma\geq\Cn_p(\Omega_\tau)^{-\frac{n-p-1}{n-p}}\int\limits_{\set{u=\tau}}\abs{\D u_\tau}^p \dd \sigma\label{eq:maininequalitymonotonicity}
\end{equation}
where $u_\tau= u/\tau $. The H\"{o}lder inequality with conjugate exponents $a=(p+1)/p$ and $b=p+1$, yields
\begin{equation}
    \capa_p(\Omega_\tau)^{\frac{p+1}{p}} \leq \left(\kern.1cm \int\limits_{\set{u=\tau}} \abs{\D u_\tau }^p \dd \sigma \right) \left( \kern.1cm\int\limits_{\set{u=\tau}} \frac{1}{\abs{\D u_\tau} } \dd \sigma \right)^{\frac{1}{p}}.
\end{equation}
Therefore, plugging it into \eqref{eq:maininequalitymonotonicity}, we get
\begin{equation}
    \abs{ \Sf^{n-1}}\Cn_p(\Omega_\tau)^{\frac{n}{n-p} }\leq \left( \frac{n-p}{p-1}\right) \otheta^{\frac{p}{n-p}} \int\limits_{\set{u=\tau}} \frac{ 1 }{\abs{\D u_\tau}} \dd \sigma
\end{equation}
Using \eqref{eq:scaling-invariant-cap} and integrating both sides we obtain 
\begin{align}
\abs{\Sf^{n-1}}\Cn_p(\Omega)^{\frac{n}{n-p}}\int_\tau^1 s^{-\frac{n(p-1)}{n-p}-1}\dd s\leq \left( \frac{n-p}{p-1}\right)\otheta^{\frac{p}{n-p}}\int_\tau^1 \int\limits_{\set{u=s}} \frac{1}{\abs{\D u}} \dd \sigma \dd s
\end{align}
that, together with the coarea formula \eqref{eq:enached_coarea} with $f= (1-\chi_{\Crit(u)}) \abs{\D u}^{-1} $ (see \cref{rmk:coarea_for_fraction}), leaves us with
\begin{align}
\frac{\abs{\Sf^{n-1}}}{n}\left(\Cn_p(\Omega_\tau)^{\frac{n}{n-p}}-\Cn_p(\Omega)^{\frac{n}{n-p}}\right)
\leq\otheta^{\frac{p}{n-p}} \abs{ \Omega_\tau \smallsetminus (\Omega\cup \Crit(u))}
\end{align}
for every $\tau \in [0,1)$. Applying the sharp iso-$p$-capacitary inequality \eqref{eq:isopcap} to the left hand side we obtain
\begin{equation}
    \AVR(g)^{\frac{p}{n-p}}\left( \abs{ \Omega_\tau}- \Cn_p(\Omega)^{\frac{n}{n-p}}\right)\leq \otheta^{\frac{p}{n-p}} \abs{ \Omega_\tau}
\end{equation}
Dividing both sides by $\abs{\Omega_\tau}$ and passing to the limit as $\tau\to 0$, we get a contradiction with $\otheta < \AVR(g)$, proving that for any $\otheta < \AVR(g)$
\begin{equation}\label{eq:Minkowski-contradicition-theta}
    \mathrm{C}_p(\Omega)^{\frac{n-p-1}{n-p}}\otheta^{\frac{1}{n-p}} < \frac{1}{\abs{\S^{n-1}}} \left(\frac{p-1}{n-p}\right)^p \int_{\partial \Omega} \abs{ \D u }^p \dd \sigma
\end{equation}
holds for every any bounded open $\Omega \subset M$  with smooth boundary.
Letting $\otheta \to \AVR(g)^{-}$ yields \eqref{eq:Minkowski-contradicition}.

To conclude observe that \cref{prop:Monotonicity_Fbeta_sequence} implies $(F_p)'(1) \leq 0$ and thus, thanks to \eqref{eq:derivative_of_F_standard}, we have 
\begin{align}
\int\limits_{\partial \Omega} \left(\frac{p-1}{n-p}\right) \abs{\D u}^p \dd \sigma\leq \int\limits_{\partial \Omega} \abs{\D u }^{p-1} \frac{\HH}{n-1}\dd \sigma.
\end{align}
By H\"{o}lder inequality with conjugate exponents $a=p/(p-1)$ and $b=p$, we get
\begin{equation}
\int\limits_{\partial \Omega} \abs{\D u}^p \dd \sigma \leq \left( \frac{n-p}{p-1}\right)^p \int\limits_{\partial \Omega} \abs{ \frac{\HH}{n-1}}^p \dd \sigma,\label{eq:LpMinkowskiHolder}
\end{equation}
that coupled with \eqref{eq:Minkowski-contradicition} concludes the proof of \eqref{eq:Lp_minkowski}.
\smallskip

If we now assume that the equality holds in \eqref{eq:Lp_minkowski}, then the two sides of \eqref{eq:LpMinkowskiHolder} are identical too. In particular, by \eqref{eq:derivative_of_F_standard}, $F'_p(1)=0$ and the rigidity statement in \cref{prop:Monotonicity_Fbeta_sequence} applies.
\endproof

In order to derive the Extended Minkowski Inequality we want to briefly recall the definition of outward minimising sets and the notion of strictly outward minimising hull in accordance to \cite{Huisken2001} and some related properties that the interested reader can find in \cite{Fogagnolo2020a}. We are denoting with $\partial^*E$ the reduced boundary of a finite perimeter set $E$. 

\begin{definition}[Outward minimising and strictly outward minimising sets] 
Let $(M,g)$ be a complete Riemannian manifold. Let $E \subset M$ be a bounded measurable set with finite perimeter. $E$ is \emph{outward minimising} if for any $F\supseteq E$ we have $\abs{\partial^* E} \leq \abs{\partial^* F}$, where by $\partial^* F$ we denote the reduced boundary of a set $F $. $E$ is \emph{strictly outward minimising} if it is outward minimising and whenever $\abs{\partial^* E}=\abs{\partial^*F}$ for some $F \supseteq E$ we have that $\abs{F \smallsetminus E} =0$.
\end{definition}
 
We can define the \emph{strictly outward minimising hull} $\Omega^*$ of an open bounded subset $\Omega$ with smooth boundary as
\begin{align}\label{eq:def_SOMH}
\Omega^* =\Int E && \begin{aligned}[c] &\text{for some bounded $E$ containing $\Omega$}\\ &\text{such that }\abs{E} = \inf_ {F \in \mathrm{SOMBE}(\Omega)} \abs{F},\end{aligned}
\end{align}
where by $\mathrm{SOMBE}(\Omega)$ we denote the family of all bounded strictly outward minimising sets containing $\Omega$ and $\Int E$ is the measure theoretic interior of $E$. As a consequence of \cite[Theorem 1.1]{Fogagnolo2020a}, if $(M, g)$ is a manifold with nonnegative Ricci curvature and Euclidean volume growth, then $\Omega^*$ as defined above is unique and it is a maximal volume solution to the problem of area minimisation with obstacle $\Omega$, that is 
\begin{equation}\label{eq:minimizing_problem}
\abs{\partial^* \Omega^*} = \inf\set{\abs{\partial^*F} \st F \text{ is bounded and } \Omega \subseteq F}.
\end{equation}
Outward minimising sets can be characterised as those satisfying
\begin{equation}\label{eq:caracterization_of_OMS}
\abs{\partial \Omega}= \abs{ \partial \Omega^*}.
\end{equation}
The relation between the strictly outward minimising hull of a bounded set with smooth boundary $\Omega$ and its $p$-capacity in the family of manifolds we are working on is resumed in the limit
\begin{align}\label{eq:limit_of_p_capacitary}
\lim_{p\to 1^+} \mathrm{C}_p(\Omega) = \frac{\abs{\partial \Omega^*}}{\abs{\S^{n-1}}}.
\end{align}
Such result is contained in the far more general \cite[Theorem 1.2]{Fogagnolo2020a}, having in mind the relation between the $p$-capacity and the normalised $p$-capacity given in \cref{def:capandnormpcap}. Letting $p \to 1^+$ in the $L^p$-Minkowski Inequality \eqref{eq:Lp_minkowski} and employing the Dominated Convergence Theorem complete the proof of the Extended Minkowski Inequality of \cref{thm:EMI_maintheorem}
\begin{equation}\label{eq:EMI_text}
\left(\frac{\abs{\partial \Omega^*}}{\abs{\S^{n-1}}}\right)^{\frac{n-2}{n-1}}\AVR(g)^{\frac{1}{n-1}}\leq \frac{1}{\abs{\S^{n-1}}}\int\limits_{\partial \Omega} \abs{ \frac{\HH}{n-1}} \dd \sigma.
\end{equation}

Outward minimising sets are mean-convex, as a simple variational argument immediately shows, and satisfy \eqref{eq:caracterization_of_OMS}. As a corollary, the Minkowski Inequality can be simplified for this particular class of subsets as in the following statement. 

\begin{corollary}[Minkowski Inequality for outward minimising sets]
\label{corollary-minimising}
Let $(M,g)$ be complete Riemannian manifold with $\Ric \geq 0$ and Euclidean Volume Growth. Let $\Omega \subseteq M$ be a bounded outward minimising subset with smooth boundary, then
\begin{align}
\label{minkonical-outward}
\left(\frac{\abs{\partial \Omega}}{\abs{\S^{n-1}}}\right)^{\frac{n-2}{n-1}} \AVR(g)^{\frac{1}{n-1}}\leq \frac{1}{\abs{\Sf^{n-1}}}\int\limits_{\partial \Omega} \frac{\HH}{n-1} \dd \sigma.
\end{align}
\end{corollary}

\subsection{Rigidity statement} 
We finally characterise the subsets $\Omega$ that saturate the inequality \eqref{eq:EMI_text}. We are getting this rigidity result evolving $\partial \Omega$ by smooth IMCF, proving that, in a outer neighbourhood of $\partial \Omega$, the manifolds is a truncated cone with the same volume ratio of $(M,g)$. The conclusion then follows from a generalisation of the Bishop-Gromov Theorem.

Going into more detail, since $\partial \Omega$ is strictly mean-convex, we can consider a sequence of sets $\Omega_t$ with $t\in [0, T)$ such that $\partial \Omega_t = F_t(\partial \Omega)$ where $F_t : \partial \Omega \to M$ satisfies
\begin{equation}
\label{imcf-def}
\frac{\!\dd}{\!\dd t} F_t (\partial \Omega) = \frac{1}{\HH_t} \nu_t,
\end{equation}
where $\nu_t$ and $\HH_t$ are respectively the outer unit normal and the mean curvature of $\partial \Omega_t$.
The conical splitting we aim to is inspired by an argument contained in
\cite[Section 8]{Huisken2001}. A first step consists in the following fundamental Lemma.

\begin{lemma} \label{lem:CMC_lemma}
Let $(M,g)$ be a complete Riemannian manifold with $\Ric\geq 0$ and $\Sigma \subseteq M$ a totally umbilical closed hypersuface such that $\Ric(\nu, \nu )=0$ where $\nu$ is the normal unit vector field to $\Sigma$. Then $\Sigma$ has \emph{constant} mean curvature.
\end{lemma}

\proof The (traced) Codazzi-Mainardi equations and the totally umbilicity yields
\begin{equation}
    \Ric_{j\nu} =\D_i \hh_{ij}- \D_j \HH= -\frac{n-2}{n-1} \D_j \HH
\end{equation}
for any $j= 1, \ldots, n-1$. Consider, at a fixed point on $\Sigma$, the vector $\eta_\lambda = \lambda\, \D^\top \HH + \nu$, with $\lambda \in \R$. Since $\Ric(\nu,\nu)=0$, we have
\begin{equation}
0 \leq \Ric(\eta_\lambda ,\eta_\lambda)=2 \Ric_{j \nu}\eta_\lambda^j\eta_\lambda^\nu +\Ric_{ij} \eta_\lambda^i \eta_\lambda^j = -2\lambda \frac{n-2}{n-1}  \abs{ \D^ \top \HH}^2 +\lambda^2\Ric_{ij} \D^i \HH\, \D^j \HH
\end{equation}
for every $\lambda \in \R$. This can happen only if $\abs{ \D^\top \HH}=0$, so that $\HH$ is constant on $\Sigma$. \endproof
The following straightforward but very important consequence of the  Bishop-Gromov monotonicity ensures in particular that if an outer neighbourhood of a bounded open set with smooth boundary $\Omega \subset M$ is isometric to a truncated cone, then the whole complement of $\Omega$ is isometric to a truncated cone based at $\partial \Omega$.

\begin{lemma}\label{lem:BG_conico}
Let $(M,g)$ be a complete noncompact Riemannian manifold with $\Ric \geq 0$. Let $K \subset M$ be a bounded open set. Suppose there exists an outer neighbourhood $A \subset M\smallsetminus K$ of $K$  such that $(A, g)$ is isometric to
\begin{equation}
    \left([\rho_0, \rho_1]\times \partial K, \dd \rho \otimes \dd \rho + \left(\frac{\rho}{\rho_0}\right)^2 g_{\partial K}\right)
\end{equation}
for $0<\rho_0 <\rho_1$. Then
\begin{equation}\label{eq:BG_pieceofcone}
    \abs{\partial K} \geq\rho_0^{n-1}\abs{\S^{n-1}}\AVR(g),
\end{equation}
and the equality holds if and only $(M \smallsetminus K, g)$ is isometric to 
\begin{equation}
    \left([\rho_0, + \infty)\times \partial K, \dd \rho \otimes \dd \rho + \left(\frac{\rho}{\rho_0}\right)^2 g_{\partial K}\right).
\end{equation}

\end{lemma}
\begin{proof}
Consider the cone $(C, \hat{g})$ given by 
\[
\left((0, \rho_1)\times \partial K, \dd \rho \otimes \dd \rho + \left(\frac{\rho}{\rho_0}\right)^2 g_{\partial K}\right),
\]
and the Riemannian manifold, with a conical singularity, obtained by gluing $(C, \hat{g})$ with $(M \smallsetminus (K \cup A), g)$ along $\{\rho = \rho_1\}$. By our assumptions, such  manifold is well-defined with nonnegative Ricci curvature outside of the tip $o$ of $C$, and coincides with $(M, g)$ in the complement of $K$.
In $C$, the geodesic distance from $o$ is given by $\rho$, and in particular, by Bishop-Gromov monotonicity,
\[
\frac{\abs{\set{\rho = r}}}{r^{n-1} \abs{\Sf^{n-1}}} \geq \AVR(g),
\]
for any $r \in (0, \rho_1)$. Since ${\abs{\set{\rho = \rho_0}}} = \abs{\partial K}$, setting $r = \rho_0$ proves \eqref{eq:BG_pieceofcone}. If equality holds, then, by the rigidity statement in Bishop-Gromov Theorem for manifolds with a conical singularity, the whole manifold we constructed is isometric to a cone, and in particular $(M \smallsetminus K, g)$ splits as claimed. This well-known, slightly enhanced version of Bishop-Gromov rigidity statement can be readily deduced from its classic proof, or seen as a very special case of its version for nonsmooth metric spaces \cite{Dephilippis2016}.
\end{proof}
We finally have at our disposal all the tools we need to work out the splitting argument leading to \cref{thm:rigidity}.

\medskip
\proof[Proof of \cref{thm:rigidity}] Suppose that some strictly outward minimising $\Omega\subset M$ with strictly mean-convex boundary satisfies
\begin{equation}\label{minkonical-equality}
\left(\frac{\abs{\partial \Omega}}{\abs{\S^{n-1}}}\right)^{\kern-.05cm\frac{n-2}{n-1}}\kern-.1cm\AVR(g)^{\frac{1}{n-1}}= \frac{1}{\abs{\S^{n-1}}}\int\limits_{\partial \Omega}  \frac{\HH}{n-1}\dd \sigma.
\end{equation}
Since $\partial \Omega $ is by assumption strictly mean-convex, we can evolve it by (smooth) IMCF $\partial \Omega_t$ defined in \eqref{imcf-def} for $t\in [0,T)$. By the \cite[Smooth Start Lemma 2.4]{Huisken2006}, up to shortening the time interval, we can assume that $\Omega_t$ is strictly outward minimising for any $t \in [0,T)$. Indeed, since $\Omega$ is strictly outward minimising, the flow coincides for a short time with the weak notion of IMCF, that exists in our setting by \cite[Theorem 1.8]{Mari2019}. The sublevel sets of the weak IMCF being strictly outward minimising is a basic and fundamental property illustrated in \cite[Minimizing Hull Property 1.4]{Huisken2001}. Consider then the function $\mathcal{Q}:[0,T) \to \R$ defined by 
\begin{equation}
\label{q-function}
 \mathcal{Q} (t) \, = \, 
|\pa\Om_t|^{-\frac{n-2}{n-1}} \!\int\limits_{\pa\Om_t} \!\HH_t \,\dd \sigma\, .
\end{equation}
A straightforward computation, direct consequence of the evolution equations for curvature flows derived for example in \cite[Theorem 3.2]{Huisken1999} shows that
\begin{equation}
\label{q-der}
\mathcal{Q}'(t) = - \abs{\partial \Omega_t}^{-\frac{n-2}{n-1}} \int\limits_{\partial \Omega_t} \frac{\abs{\mathring{\hh}_t}^2 + \ric(\nu_t, \nu_t)}{\HH_t} \dd \sigma\leq 0 \,,
\end{equation} 
where by $\mathring{\hh}_t$ we denote the trace-free part of the second fundamental form $\hh_t$ of $\partial \Omega_t$. On the other hand, the strict inequality for some $t \in [0,T)$ would result in a contradiction to the Minkowski Inequality. Thus $\mathcal{Q}'(t)$ vanishes for any $t \in [0,T)$ and, in particular $\partial \Omega_t$ satisfies \eqref{minkonical-equality} for any $t \in [0,T)$. Hence, $\partial \Omega_t$ is totally  umbilical and satisfies $\Ric(\nu_t,\nu_t) = 0$ in for every $t \in [0,T)$. By \cref{lem:CMC_lemma} $\partial \Omega_t$ has constant mean curvature for every $t\in [0,T)$. 

On $\set{0\leq t <T}$ the function $w$ 
solution to the weak level set formulation of the IMCF, that in our smooth case just means $\{w = t\} = \partial \Omega_t$,
satisfies the relation
\begin{equation}
\label{weak-imcf-eq}
\HH_{t} = \mathrm{div}\left(\frac{\D w}{\abs{\D w}}\right)(x_t) = \abs{\D w}(x_t)
\end{equation}
at any $x_t \in \partial \Omega_t$.
Hence, since $\HH_t>0$, a well-known extension of the Gauss' Lemma yields
\begin{equation}\label{splitting}
g = \frac{\!\dd w \otimes \! \dd w}{\abs{\D w}^2} + g^{}_{\partial \Omega_{\{w = t\}}} = \frac{\!\dd t \otimes \! \dd t}{\HH_t^2} +g^{}_{\partial \Omega_t}.
\end{equation}
The evolution equation (see \cite[Theorem 3.2 (i)]{Huisken1999}) satisfied by $g^{}_{\partial \Omega_t}$ is 
\begin{equation}
\label{evolution-metric}
\frac{\partial}{\partial t} g^{}_{\partial \Omega_t} = 2\frac{\hh_t}{\HH_t} =\frac{2}{n-1} g^{}_{\partial \Omega_t},
\end{equation}
where the last identity is due to the total umbilicity of $\partial \Omega_t$. Integrating such equation we deduce 
\begin{equation}
\label{formula-metric}
g^{}_{\partial \Omega_t} = \mathrm{e}^{\frac{2t}{(n-1)}}g^{}_{\partial \Omega},
\end{equation}
On the other hand, the evolution equation for the mean curvature along the IMCF (see \cite[Theorem 3.2 (v)]{Huisken1999}) declaims
\begin{equation}
\label{evolution-mean}
\frac{\partial}{\partial t} \HH_t = - \Delta_{\partial \Omega_t}\left(\frac{1}{\HH_t}\right) - \frac{1}{\HH_t}\left[\abs{\hh_t}^2 + \ric (\nu_t, \nu_t)\right] = -\frac{\HH_t}{n-1}
\end{equation}
where the last identity is due to the fact that $\partial \Omega_t$ is totally umbilical, $\Ric(\nu_t,\nu_t)=0$ and the mean curvature $\HH_t$ of $\partial \Omega_t$ depends only on $t$. Integrating it we obtain that
\begin{equation}
\label{formula-mean-curvature}
    \HH_t = \ee^{-\frac{t}{n-1}}\HH_0
\end{equation}
where $\HH_0$ is the mean curvature of $\partial \Omega$.

Plugging \eqref{formula-metric} and \eqref{formula-mean-curvature} into \eqref{splitting}, we deduce that $(\set{0\leq t < T},g)$ is isometric to
\begin{equation}
    \left([0,T) \times \partial \Omega,  \, \ee^{\frac{2t}{n-1}}\frac{\dd t\otimes \dd t}{\HH_0^2} + \ee^{\frac{2t}{n-1}} g_{\partial \Omega}\right).
\end{equation}
Performing the change of variables
\begin{equation}
    {\rho} = \frac{(n-1)}{\HH_0}\mathrm{e}^{\frac{t}{(n-1)}},
\end{equation}
the metric can be written as
\begin{align}
\left([\rho_0, \rho(T)) \times \partial \Omega, \dd \rho\otimes \dd \rho + \left(\frac{\rho}{\rho_0}\right)^2g_{\partial \Omega}\right), && \text{where } \rho _0 = \frac{(n-1)}{\HH_0} .
\end{align}
On the other hand, since by assumption $\partial \Omega$ saturates the Minkowski Inequality, that is \eqref{minkonical-equality} holds, we immediately get 
\[
\rho _0 = \left(\frac{\abs{\partial \Omega}}{\AVR(g)\abs{\S^{n-1}}}\right)^{\frac{1}{n-1}},
\]
and we conclude by the rigidity statement in \cref{lem:BG_conico} that the whole $M \smallsetminus \Omega$ is isometric to a truncated cone.
\endproof 
In the following Remark we briefly discuss how the assumptions for the rigidity can be relaxed in small dimensions.
\begin{remark}
In dimension $3\leq n \leq 7$, an open bounded subset $\Omega$ with smooth strictly mean-convex boundary satisfying
\begin{align}
\left(\frac{\abs{\partial \Omega^*}}{\abs{\S^{n-1}}}\right)^{\!\!\frac{n-2}{n-1}} \AVR(g)^{\frac{1}{n-1}} = \frac{1}{\abs{\Sf^{n-1}}} \int\limits_{\partial \Omega} {\frac{\HH}{n-1}} \dd \sigma.
\end{align}
is \emph{a priori} strictly outward minimising, and thus, in this case, such assumption can be dropped. Indeed, by approximating $\Omega$ via mean curvature flow with smooth strictly outward minimizing domains, as described in \cite[Lemma 5.6]{Huisken2001}, we deduce that \eqref{eq:EMI_text} holds also for $\CC^{1,1}$-hypersurfaces. In particular, the Minkowski Inequality holds also for the strictly outward minimising hull of $\Omega$ (see the regularity results recalled in \cite[Regularity Theorem 1.3]{Huisken2001} and \cite[Theorem 2.18]{Fogagnolo2020a} and in the references therein) for every $\Omega$ with smooth boundary, provided the dimensional bound holds. We can then argue by contradiction. Suppose that $\Omega^*$ does not coincide with $\Omega$, then
\begin{align}
\left(\frac{\abs{\partial \Omega^*}}{\abs{\S^{n-1}}}\right)^{\!\!\frac{n-2}{n-1}} \AVR(g)^{\frac{1}{n-1}} = \frac{1}{\abs{\Sf^{n-1}}} \int\limits_{\partial \Omega} {\frac{\HH}{n-1}} \dd \sigma > \frac{1}{\abs{\Sf^{n-1}}} \int\limits_{\partial \Omega^*} {\frac{\HH}{n-1}} \dd \sigma
\end{align} 
where the last inequality is due to the fact that $\HH=0$ on $\partial \Omega^* \smallsetminus \partial \Omega$. But this contradicts the Minkowski Inequality for $\Omega^*$, hence $\Omega= \Omega^*$. 
\end{remark}

\subsection{A pinching condition and a sphere theorem}
In this subsection, we exploit the monotonicity of the function $t \mapsto F_p^\infty(t)$ defined in \eqref{eq:monotone_function_infty} to prove a couple of rigidity statements involving a pinching condition on the mean curvature of $\partial \Omega$ and an {\em a priori} bound on the gradient of the $p$-capacitary potential associated to $\Omega$. These results without any convexity assumption are new also in $\R^n$, and they constitute the complete nonlinear generalisation of \cite[Corollary 1.4 and 1.9]{Borghini2019}. For convex subsets of the Euclidean space they are the content of \cite[Corollary 2.16 and 2.17]{Fogagnolo2019}.

\begin{theorem}\label{thm:p-pinching}
Let $(M,g)$ be a complete Riemannian manifold with $\Ric\geq 0$ and Euclidean Volume Growth. If there exists a open bounded subset $\Omega \subseteq M$ with smooth boundary satisfying
\begin{equation}\label{eq:p-pinching_condition}
-\left[\frac{\AVR(g)}{\mathrm{C}_p(\Omega)}\right]^{\frac{1}{n-p}}\leq \frac{\HH}{n-1}\leq\left[\frac{\AVR(g)}{\mathrm{C}_p(\Omega)}\right]^{\frac{1}{n-p}}
\end{equation}
on every point of $\partial \Omega$, then $(M \smallsetminus \Omega, g)$ is isometric to
\begin{align}
\left([\rho_0,+\infty) \times \partial \Omega, \dd \rho \otimes \dd \rho  + \left(\frac{\rho }{\rho_0}\right)^2 g^{}_{\partial \Omega}\right),&& \text{where } \rho _0 = \left(\frac{\abs{\partial \Omega}}{\AVR(g)\abs{\S^{n-1}}}\right)^{\frac{1}{n-1}}.
\end{align}
In this case $\partial \Omega$ is a connected totally umbilical hypersurface with constant mean curvature in $(M\smallsetminus \Omega, g)$.
\end{theorem}

\proof We can argue by contradiction as in \cref{thm:Lp_mink} to prove that
\begin{equation}
    \left(\frac{n-p}{p-1}\right) \left[\frac{\AVR(g)}{\mathrm{C}_p(\Omega)}\right]^{\frac{1}{n-p}}\leq \sup_{\partial \Omega} \abs{ \D u}.
\end{equation}
Indeed, we can follow the same lines replacing the consequence of the monotonicity of $F_p$ with the corresponding of $F_p^\infty$, that thanks to \eqref{eq:scaling-invariant-cap} can be rewritten as
\begin{equation} \label{eq:monotonicity_infty_geometric}
    F^{\infty}_{p}(t)= t^{\frac{n-1}{n-p}} \sup_{\set{u=1/t}} \abs{ \D u} = \left( \frac{\Cn_p(\Omega_t)}{\Cn_p(\Omega)}\right)^{\frac{1}{n-p}} \sup_{\set{u_t=1}} \abs{\D u_t}
\end{equation}
where $u_t=tu$ and $\Omega_t= \set{u>1/t}\cup \Omega$. Accordingly, we employ the H\"{o}lder inequality with conjugate exponents $a=+\infty$ and $b=1$, that is
\begin{equation}
 \capa_p(\Omega_t)^{\frac{1}{p}} \leq\sup_{\set{u=1/t}} \abs{ \D u_t} \left( \int\limits_{\set{u=1/t}} \frac{1}{\abs{\D u_t}} \dd \sigma  \right)^{\frac{1}{p}}.
\end{equation}
In the end, by \eqref{prop:Monotonicity_Finfty_sequence} we get
\begin{equation}
\sup_{\partial \Omega} \abs{ \D u} \leq  \frac{(n-p)}{(p-1)(n-1)} \sup_{\partial \Omega} \abs{\HH}
\end{equation}
and the equality holds if and only if $(M\smallsetminus \Omega, g)$ splits as in the statement. Condition \eqref{eq:p-pinching_condition} easily implies the equality. \endproof

The above result is a rigidity theorem under a pinching condition on the mean curvature of $\partial \Omega$ with respect to its $p$-capacity. From the proof above we can also get that
\begin{equation}\label{eq:cone_condtion}
\frac{1}{p-1} \left[\frac{\AVR(g)}{\mathrm{C}_p(\Omega)}\right]^{\frac{1}{n-p}}\leq \sup_{\partial \Omega} \abs{ \frac{\D u}{n-p}}
\end{equation}
and the equality is satisfied  only on metric cones. The previous inequality gives a lower bound on the gradient of $u$ on $\partial \Omega$ in terms of the $p$-capacity of $\Omega$ that, when attained, forces $(M,g)$ to be (isometric to) $\R^n$ with $\Omega$ a Euclidean ball.

\begin{theorem}
\label{thm:gradient pinching}
Let $(M,g)$ be a complete Riemannian manifold with $\Ric\geq 0$ curvature and Euclidean Volume Growth. Let $\Omega \subseteq M$ be an open bounded subset with smooth boundary, $u$ the $p$-capacitary potential associated to $\Omega$ and assume that
\begin{equation}\label{eq:condition_euclidean}
\sup_{\partial \Omega} \abs{ \frac{\D u}{n-p}} \leq \frac{1}{p-1} \AVR(g)^{\frac{1}{p-1}}\left(\frac{\abs{\S^{n-1}}}{\abs{\partial \Omega}}\right)^{\frac{1}{n-1}}.
\end{equation}
Then $(M,g)$ is isometric to $\R^n$ with the Euclidean metric and $\Omega$ is ball.
\end{theorem}

\proof Under the assumption \eqref{eq:condition_euclidean}, we get
\begin{equation}
\mathrm{C}_p(\Omega) =\left(\frac{p-1}{n-p}\right)^{p-1} \frac{1}{\abs{\S^{n-1}}} \int\limits_{\partial \Omega} \abs{ \D u }^{p-1} \dd \sigma \leq \AVR(g) \left(\frac{\abs{\S^{n-1}}}{\abs{ \partial \Omega}}\right)^{-\frac{n-p}{n-1}}
\end{equation}
that yields
\begin{equation} 
\left(\frac{\abs{\S^{n-1}}}{\abs{ \partial \Omega}}\right)^{\frac{n-p}{n-1}}\leq \frac{\AVR(g)}{\mathrm{C}_p(\Omega)}\leq (p-1)^{n-p} \sup_{\partial \Omega} \abs{ \frac{\D u}{n-p}}^{n-p} \leq \AVR(g)^{\frac{n-p}{p-1}}\left( \frac{\abs{\S^{n-1}}}{\abs{ \partial \Omega}}\right)^{\frac{n-p}{n-1}}\label{eq:condition_on_Omega}
\end{equation}
where we used \eqref{eq:cone_condtion} together with the condition \eqref{eq:condition_euclidean}. Thus, we obtain that $\AVR(g)=1$, and hence, by Bishop-Gromov Theorem, that $(M,g)$ is isometric to $\R^n$ with the Euclidean metric. Since all inequalities in \eqref{eq:condition_on_Omega} becomes equalities, by the second one we can apply the rigidity statement in \cref{prop:Monotonicity_Finfty_sequence} which ensures that $\partial \Omega$ is a compact connected and totally umbilical hypersurface of $\R^n$, that is, $\Omega$ is a ball. \endproof

\bigskip
 \appto\bibfont{\setlength{\emergencystretch}{1em}}
 \printbibliography

\end{document}